\documentclass[12pt]{article}

\RequirePackage[OT1]{fontenc}
\RequirePackage{amsthm,amsmath,natbib,amssymb,latexsym,bbm}

\newfont{\eufb}{eufb10 scaled\magstep1}
\newfont{\eufm}{eufm10 scaled\magstep1}
\newfont{\eusb}{eusb10 scaled\magstep1}

\newcommand{\DS}{\displaystyle}

\newcommand{\SC}{\scriptstyle}

\newcommand{\NN}{\mathbb{N}}
\newcommand{\RR}{\mathbb{R}}
\newcommand{\ZZ}{\mathbb{Z}}

\newcommand{\EE}{\operatorname{\mathsf{E}}}
\newcommand{\PP}{\operatorname{\mathsf{P}}}
\newcommand{\bff}{\mathbf{f}}
\newcommand{\bg}{\mathbf{g}}
\newcommand{\bh}{\mathbf{h}}
\newcommand{\tDelta}{\widetilde{\Delta}}
\newcommand{\tD}{\widetilde{D}}
\newcommand{\tG}{\Tilde{G}}
\newcommand{\ttG}{\Tilde{\tG}}
\newcommand{\tI}{\Tilde{I}}
\newcommand{\tJ}{\Tilde{J}}
\newcommand{\ttI}{\Tilde{\tI}}
\newcommand{\tM}{\widetilde{M}}
\newcommand{\tR}{\widetilde{R}}
\newcommand{\tx}{\widetilde{x}}

\newcommand{\dd}{\mathrm{d}}
\newcommand{\ee}{\mathrm{e}}
\newcommand{\ii}{\mathrm{i}}

\newcommand{\vare}{\varepsilon}

\renewcommand{\Re}{\operatorname{\mbox{\eufm Re}}}
\renewcommand{\Im}{\operatorname{\mbox{\eufm Im}}}

\numberwithin{equation}{section}
\theoremstyle{plain}
\newtheorem{thm}{Theorem}[section]
\newtheorem{pro}{Proposition}[section]
\newtheorem{lem}{Lemma}[section]
\theoremstyle{remark}

\begin{document}

\title{The accuracy of merging approximation in generalized St.~Petersburg
       games \\
Running title: Accuracy of merging approximation}

\author{Gyula Pap \\
        Faculty of Informatics,
        University of Debrecen, \\
        Pf.12,
        H--4010 Debrecen,
        Hungary}

\date{April 14, 2009}

\maketitle

\begin{abstract}
Merging asymptotic expansions of arbitrary length are established for the
 distribution functions and for the probabilities of suitably centered and
 normalized cumulative winnings in a full sequence of generalized
 St.~Petersburg games, extending the short expansions due to Cs\"org\H{o}, S.,
 Merging asymptotic expansions in generalized St.~Petersburg games,
 \textit{Acta Sci.\ Math.\ (Szeged)} \textbf{73} 297--331, 2007.
These expansions are given in terms of suitably chosen members from the classes
 of subsequential semistable infinitely divisible asymptotic distribution
 functions and certain derivatives of these functions.
The length of the expansion depends upon the tail parameter. 
Both uniform and nonuniform bounds are presented.
\vspace*{3mm}

\vspace*{3mm}

\end{abstract}

\section{Introduction}
\label{intro}

For \ $ p \in (0,1) $ \ and \ $ \alpha > 0 $, \ consider a generalized
 St.~Petersburg game in which the gain \ $ X $ \ of a gambler is such that
 \ $ \PP \left\{ X = r^{k / \alpha} \right\} = q^{k-1} p $ \ for
 \ $ k \in \NN := \{ 1, 2, \dots\} $, \ where \ $ q := 1 - p $ \ and
 \ $ r := 1 / q $.
\ We refer to Cs\"org\H{o} \cite{Cso_02} for a brief introduction to both the
 history and the mathematics of this game.
Let \ $ X_1 $, $ X_2 $, \dots \ be the gambler's gains in a sequence of
 independent repetitions of the game, with \ $ S_n := X_1 + \cdots + X_n $
 \ standing for the total winnings in \ $ n $ \  games, \ $ n \in \NN $.
\ The main purpose of the present paper is to study the asymptotic behaviour of
 the distributions of \ $ S_n $.
\ Note that \ $ \EE( X^\alpha ) = \infty $ \ and \ $ \EE( X^\beta ) < \infty $
 \ for all \ $\beta \in ( 0, \alpha) $, \ namely,
 \[
   \EE( X^\beta ) = \frac{p} { q^{\beta / \alpha} - q } =: \mu_\beta^{ \alpha, p } .
 \] 
Remark that even for \ $ \beta > \alpha $, \ a virtual moment
 \ $ \mu_\beta^{ \alpha, p } < 0 $ \ may be defined by this formula.
In the finite-variance case \ $ \alpha > 2 $, \ the central limit theorem is
 valid, namely,
 \begin{equation} \label{CLT}
  \lim_{n \to \infty}
   \PP \bigg\{ \frac{S_n - \mu_1^{ \alpha, p } n } { a_n^{ \alpha, p } } \leq x
       \bigg\}
   = \Phi(x) \qquad \text{for all \ $ x \in \RR $,}
 \end{equation}
 where \ $ \RR $ \ is the real line,
 \ $ a_n^{ \alpha, p } := \sigma_{\alpha, p} \sqrt{n}$ \ with
 \ $\sigma_{\alpha, p} := \sqrt{ \mu_2^{ \alpha, p } - [\mu_1^{ \alpha, p }]^2 } $, \ and
 \ $ \Phi $ \ is the standard normal distribution function.
In the infinite-variance case \ $ \alpha = 2 $ \ the underlying distribution is
 still in the domain of attraction of the normal law and \eqref{CLT} holds with
 the choice \ $ a_n^{ \alpha, p } := \sqrt{ p r \, n \log_r n } $, \ see
 Cs\"org\H{o} \cite{Cso_02}.
By standard results for normal approximation, one can derive rate of
 convergence in the central limit theorem \eqref{CLT}.
Namely, there exist constants \ $ C(\alpha,\beta,p) > 0 $ \ and
 \ $ C(\alpha,p) > 0 $ \ such that for all \ $n \in \NN$ \ and \ $x \in \RR$,
 \[
   \left| \PP \bigg\{ \frac{S_n - \mu_1^{ \alpha, p } n }
                           { a_n^{ \alpha, p } } \leq x \bigg\}
          - \Phi(x) \right|
   \leq \begin{cases}
         \frac{C(\alpha,\beta,p)}{1 + |x|^\beta} \frac{1}{\sqrt{n}} ,
          & \text{if \ $ 3 \leq \beta < \alpha $,} \\[2mm]
         \frac{C(3,p)}{1 + x^2} \frac{\log_r n}{\sqrt{n}} ,
          & \text{if \ $ \alpha = 3 $,} \\[2mm]
         \frac{C(\alpha,p)}{1 + x^2} \frac{1}{n^{(\alpha - 2) / 2}} ,
          & \text{if \ $ 2 < \alpha < 3 $,} \\[2mm]
         {\SC C(2,p)} \frac{1}{\log_r n} , & \text{if \ $ \alpha = 2 $,}
        \end{cases}
 \] 
 see Osipov \cite{Osi_67} and Petrov \cite[Theorem V.13]{Pet_75} for
 \ $\alpha > 3$, \ Hall \cite[Theorem 4.1]{Hall_82} for \ $2 < \alpha \leq 3$,
 \ Hall \cite[Corollary 1]{Hall_83} for \ $\alpha = 2$, \ and 
 Cs\"org\H{o} \cite{Cso_02}.
Moreover, if \ $ 2 < \beta < \alpha \leq 3$ \ then there exists a bounded
 decreasing function \ $ \vare_\beta^{\alpha, p} : [1, \infty) \to [0, \infty) $
 \ with \ $\lim\limits_{u \to \infty} \vare_\beta^{\alpha, p}(u) = 0$ \ such that for
 all \ $n \in \NN$ \ and \ $x \in \RR$,
 \[
   \left| \PP \bigg\{ \frac{S_n - \mu_1^{ \alpha, p } n }
                           { \sigma_{ \alpha, p } \sqrt{n} } \leq x \bigg\}
          - \Phi(x) \right|
   \leq \frac{\vare_\beta^{\alpha, p}\big( \sqrt{n} ( 1 + |x| ) \big)}
             { ( 1 + |x|^\beta ) \, n^{(\beta - 2) / 2} } , 
 \]
 see Osipov and Petrov \cite{Osi_Pet_67}.
The standard results for asymptotic expansion in the central limit theorem
 \eqref{CLT} (see Petrov \cite[Theorems VI.1--5]{Pet_75}) are not
 applicable, since the characteristic function of the gain \ $X$ \ in one game
 \begin{equation} \label{f}
  \bff_{\alpha,p}(t)
  :=\EE \big( \ee^{\ii t X} \big)
  = \sum_{k = 1}^\infty 
     \ee^{ \ii t r^{k / \alpha}} q^{k - 1} p, \qquad t \in \RR ,
 \end{equation}
 does not satisfy Cram\'er's continuity condition
 \ $\limsup\limits_{t\to\infty} |\bff_{\alpha,p}(t)| < 1$.
\ Indeed, \ $\bff_{\alpha,p}$ \ is an almost periodic function (see, e.g.,
 Katznelson \cite[VI.5]{Katz_68}), hence for all \ $\vare > 0$, \ there exists
 \ $\Lambda > 0$ \ such that for every \ $T > 0$, \ there exists
 \ $\tau \in (T, T + \Lambda)$ \ such that
 \ $\sup\limits_{t \in \RR} | \bff_{\alpha,p}(t) - \bff_{\alpha,p}(t + \tau) | < \vare$.
\ Particularly, \ $| \bff_{\alpha,p}(0) - \bff_{\alpha,p}(\tau) | < \vare$,
 \ consequently, \ $\limsup\limits_{t\to\infty} |\bff_{\alpha,p}(t)| = 1$.  
\ However, for \ $\alpha > 3$, \ in the non-lattice case
 \ $r^{1 / \alpha} \notin \NN$, \ a short expansion is valid, namely, 
 \[
   \sup_{x \in \RR} \,
    ( 1 + |x|^3 )
    \left| \PP \bigg\{ \frac{S_n - \mu_1^{ \alpha, p } n }
                            { \sigma_{ \alpha, p } \sqrt{n} } \leq x \bigg\}
           - \Phi(x) - \frac{ Q_1^{\alpha, p}(x) }{ \sqrt{n}} \right|
   = o\left( \frac{1}{\sqrt{n}} \right) ,
 \]
 as \ $n \to \infty$ \ with
 \[
   Q_1^{\alpha, p}(x)
   := \frac{ \mu_3^{\alpha, p}}{ 6 \sqrt{2 \pi}
      [ \sigma_{\alpha, p} ]^3} ( 1 - x^2 ) \, \ee^{ - x^2 / 2} , \qquad x\in\RR ,
 \]
 see Bikelis \cite{Bik_66}.
In the lattice case \ $r^{1 / \alpha} \in \NN$, \ a longer expansion with the
 usual remainder term can be obtained by adding extra terms to the expansion,
 see Osipov \cite{Osi_72} and Petrov \cite[Theorem VI.6]{Pet_75}.

In the finite-variance lattice case \ $\alpha > 2$ \ and
 \ $r^{1 / \alpha} \in \NN$, \ the local limit theorem is also valid, namely, 
 \begin{equation} \label{LLT}
   \lim_{n \to \infty}
    \left[ \frac{\sqrt{n}}{r^{ 1 / \alpha}} \PP \{ S_n = s \}
           - \Phi_{\alpha, p}' \bigg( \frac{s - \mu_1^{\alpha, p} n}{\sqrt{n}} \bigg)
    \right]
   = 0
 \end{equation}
 for all \ $s \in r^{ 1 / \alpha} \NN$, \ where
 \ $\Phi_{\alpha, p}(x) := \Phi(x / \sigma_{\alpha, p})$, \ $x \in \RR$.
\ For \ $\alpha > 3$, \ the rate of convergence is
 \[
   \sup_{s \in r^{ 1 / \alpha} \NN}
    \left| \frac{\sqrt{n}}{r^{ 1 / \alpha}} \PP \{ S_n = s \}
           - \Phi_{\alpha, p}' \left( \frac{s - \mu_1^{\alpha, p} n}{\sqrt{n}} \right)
    \right|
   = O\left( \frac{1}{\sqrt{n}} \right) ,
 \]
 see Petrov \cite[Theorem VII.6]{Pet_75}.
For \ $2 < \alpha \leq 3$, \ a nonuniform bound
 \begin{multline*}
  \sup_{s \in r^{ 1 / \alpha} \NN}
   \Bigg( 1 + \bigg| \frac{s - \mu_1^{\alpha, p} n}{\sqrt{n}} \bigg|^2 \Bigg)
   \left| \frac{\sqrt{n}}{r^{ 1 / \alpha}} \PP \{ S_n = s \}
          - \Phi_{\alpha, p}' \bigg( \frac{s - \mu_1^{\alpha, p} n}{\sqrt{n}} \bigg)
   \right| \\
  = \begin{cases}
     O\left( \frac{\log_r n}{\sqrt{n}} \right)
      & \text{if \ $\alpha = 3$,} \\[2mm]
     O\left( \frac{1}{n^{(\alpha - 2) / 2}} \right)
      & \text{if \ $2 < \alpha < 3$,}
    \end{cases}
 \end{multline*}
 is available, see Hall \cite[Theorems 5.4 and 5.5]{Hall_82}.
There is also an expansion in the local limit theorem \eqref{LLT}, namely, for
 \ $3 \leq \beta < \alpha$,
 \[
   \sup_{s \in r^{ 1 / \alpha} \NN}
    \Bigg( 1 + \bigg| \frac{s - \mu_1^{\alpha, p} n}
                           {\sqrt{n}} \bigg|^{\lfloor \beta \rfloor} \Bigg)
    \left| \frac{\sqrt{n}}{r^{ 1 / \alpha}} \PP \{ S_n = s \}
           - \big( U_{\beta,n}^{\alpha, p} \big)'
             \left( \frac{s - \mu_1^{\alpha, p} n}{\sqrt{n}} \right)
    \right|
   = o\left( \frac{1}{n^{ ( \lfloor \beta \rfloor - 2 ) / 2}} \right) ,
 \]
 where \ $\lfloor y \rfloor := \max \{ k\in \ZZ : k \leq y \}$ \ denotes the
 lower integer part of a number \ $y \in \RR$, \ where
 \ $ \ZZ := \{ 0, \pm 1, \pm 2, \dots \} $, \ and
 \[
   U_{\beta, n}^{\alpha, p}(x)
   := \Phi(x)
      + \sum_{k = 1}^{ \lfloor \beta \rfloor - 2} \frac{Q_k^{\alpha, p}(x)}{n^{k / 2}}
 \]
 is the usual approximating function in the asymptotic expansion in the central
 limit theorem, where the functions \ $Q_k^{\alpha, p}$, \ $k\in\NN$, \ are
 defined in an appropriate way, see Petrov \cite[Theorem VII.16]{Pet_75}.

From now on, we assume \ $ \alpha \in (0, 2) $. 
\ Then the distribution function of \ $ X $ \ is not slowly varying at
 infinity, hence it is not in the domain of attraction of any (stable)
 distribution, that is, \ $ S_n $ \ cannot be centered and normalized to have a
 proper limit distribution, see Cs\"org\H{o} \cite{Cso_02}.
However, the sequence \ $ ( S_n - c_n^{\alpha, p} ) / n^{ 1 / \alpha} $ \ has proper
 limit distributions along some subsequences of \ $ \NN $ \ with centering
 sequence
 \begin{equation} \label{c}
   c_n^{\alpha, p} := \begin{cases}
                    \frac{ p n }{ q^{1 / \alpha} - q } = \mu_1^{ \alpha, p } n ,
                     & \text{if \ $ \alpha \neq 1 $,} \\[2mm]
                    p r \, n \log_r n , & \text{if \ $ \alpha = 1 $.} 
                   \end{cases}
 \end{equation}
The appropriate subsequences are regulated by the position parameter
 \begin{equation} \label{gamma_n}
   \gamma_n := \frac{n}{r^{\lceil \log_r n \rceil}} \in ( q, 1 ] ,
 \end{equation}
 which describes the location of \ $n = \gamma_n r^{\lceil \log_r n \rceil}$ \ between
 two consecutive powers of \ $r = 1 / q$, \ where
 \ $ \lceil y \rceil := \min \{ k \in \ZZ : k \geq y\} $ \ denotes the upper
 integer part of a number \ $ y \in \RR $.
\ The possible proper limit distributions are semistable infinitely divisible
 distributions \ $ \{ G_{\alpha, p, \gamma} : q < \gamma \leq 1 \} $ \ with exponent
 \ $\alpha$, \ where \ $ G_{\alpha, p, \gamma} $ \ can be given by its
 characteristic function
 \begin{equation} \label{g_apg}
   \bg_{\alpha, p, \gamma}(t)
   = \int_{-\infty}^\infty \ee^{\ii t x} \, \dd G_{\alpha, p, \gamma}(x)
   = \ee^{ y_\gamma^{\alpha, p}(t) } , \qquad
   t \in \RR ,
 \end{equation}
 where
 \[
   y_\gamma^{\alpha, p}(t)
   := \begin{cases}
       \sum\limits_{k = - \infty}^\infty
        \left( \exp \left\{ \frac{\ii t r^{k / \alpha}}{\gamma^{1 / \alpha}} \right\}
               -1 \right)
        \frac{p \gamma}{q r^k},
        & \hspace*{-2mm}\text{if \ $\alpha \in (0, 1)$,} \\[3mm]
       \ii t p r \log_r \! \frac{1}{\gamma}
       + \! \sum\limits_{k = - \infty}^0 \! \!
        \left( \exp \left\{ \frac{\ii t r^{k / \alpha}}{\gamma^{1 / \alpha}} \right\}
               -1 - \frac{\ii t r^{k / \alpha}}{\gamma^{1 / \alpha}} \right) \!
        \frac{p \gamma}{q r^k} \\[3mm]
       \phantom{\ii t p r \log_r \frac{1}{\gamma}}
       + \sum\limits_{k = 1}^\infty
        \left( \exp \left\{ \frac{\ii t r^{k / \alpha}}{\gamma^{1 / \alpha}} \right\}
               -1 \right)
        \frac{p \gamma}{q r^k} ,
        & \hspace*{-2mm}\text{if \ $ \alpha = 1 $,} \\[3mm]
       \sum\limits_{k = - \infty}^\infty
        \left( \exp \left\{ \frac{\ii t r^{k / \alpha}}{\gamma^{1 / \alpha}} \right\}
               -1 - \frac{\ii t r^{k / \alpha}}{\gamma^{1 / \alpha}} \right)
        \frac{p \gamma}{q r^k},
        & \hspace*{-2mm}\text{if \ $\alpha \in (1, 2)$,}
      \end{cases}
 \]
 for all \ $ t \in \RR $, \ see Cs\"org\H{o} \cite{Cso_02} and \cite{Cso_05}.
It can be shown that for any subsequence \ $ \{ n_k \}_{k=1}^\infty $ \ of
 \ $ \NN $, \ the sequence \ $ ( S_{n_k} - c_{n_k}^{\alpha, p} ) / n_k^{ 1 / \alpha} $
 \ converges weakly as \ $ k \to \infty $ \ if and only if \ there exists 
 \ $ \gamma \in ( q, 1 ]$ \ such that
 \ $ \lceil \log_r n_k \rceil - \log_r n_k \to \log_r \frac{1}{\gamma}$
 \ $\pmod 1$ \ as \ $ k \to \infty $, \ and in this case 
 \[
   \lim_{k \to \infty}
   \PP \bigg\{ \frac{S_{n_k} - c_{n_k}^{ \alpha, p } } { n_k^{ 1 / \alpha } } \leq x
       \bigg\}
   = G_{\alpha, p, \gamma}(x) \qquad \text{for all \ $ x \in \RR $,}
 \]
 see Cs\"org\H{o} \cite{Cso_07b} and Cs\"org\H{o} and Dodunekova
 \cite{Cso_Dod_91}.
It should be noted that the special classical case \ $ \alpha = 1 $,
 \ $ p = 1 / 2$ \ and \ $ n_k = 2^k$ \ has been described by Martin-L\"of
 \cite{Mar_85}. 
The above result suggests a merging approximation of the distribution of 
 \ $ ( S_n - c_n^{\alpha, p} ) / n^{ 1 / \alpha} $ \ by \ $ G_{\alpha, p, \gamma_n} $.
\ It turns out that the rates of merge is reasonably fast.
Cs\"org\H{o} \cite{Cso_02} proved that
 \[
   \sup_{ x \in \RR }
    \Bigg| \PP \bigg\{ \frac{S_n - c_n^{ \alpha, p } } { n^{ 1 / \alpha } } \leq x
               \bigg\} - G_{\alpha, p, \gamma_n}(x) \Bigg|
   = \begin{cases}
      O \big( \frac{1}{n} \big) , & \text{if \ $\alpha \in (0, 1)$,} \\[2mm]
      O \Big( \frac{ [ \log_r n ]^2}{n} \Big) ,
       & \text{if \ $ \alpha = 1 $,} \\[2mm]
      O \Big( \frac{1}{n^{(2-\alpha)/\alpha } } \Big) ,
       & \text{if \ $\alpha \in (1, 2)$.}
     \end{cases}
 \]
In fact, the constants in these bounds are given explicitely by Cs\"org\H{o}
 \cite{Cso_02}.
For \ $ \alpha \in (0,1]$, \ Cs\"org\H{o} \cite{Cso_07b} derived short
 asymptotic expansions in the merging approximations, namely, for
 \ $\alpha \in (0,1)$,
 \[
  \sup_{ x \in \RR }
   \Bigg| \PP \bigg\{ \frac{S_n - \mu_1^{ \alpha, p } n } { n^{ 1 / \alpha } } \leq x
              \bigg\} - G_{\alpha, p, \gamma_n}(x)
                       + \frac{1}{2n} G_{\alpha, p, \gamma_n}^{(0,2)}(x)
   \Bigg|
  = \begin{cases}
     O \big( \frac{1}{n^2} \big) ,
      & \text{if \ $ 0 < \alpha \leq \frac{1}{2} $,} \\[2mm]
     O \Big( \frac{1}{ n^{ 1 / \alpha } } \Big) ,
      & \text{if \ $ \frac{1}{2} < \alpha < 1 $,}
    \end{cases}
 \]
 for \ $\alpha = 1$,
 \begin{align*}
   \sup_{ x \in \RR }
    \Bigg| & \PP \bigg\{ \frac{S_n - p r \, n \log_r n}{n} \leq x \bigg\}
             \\[2mm]
           & - G_{1, p, \gamma_n}(x)
             - \frac{p r \log_r n}{n} G_{1, p, \gamma_n}^{(1,1)}(x)
             + \frac{p^2 r^2 [\log_r n]^2}{2 n} G_{1, p, \gamma_n}^{(2,0)}(x) \Bigg|
             = O \bigg( \frac{1}{n} \bigg),
 \end{align*}
 and for \ $ \frac{1}{2} < \alpha < 1 $, \ in the non-lattice case
 \ $r^{1/\alpha} \notin \NN$, 
 \[
  \sup_{ x \in \RR }
   \Bigg| \PP \bigg\{ \frac{S_n - \mu_1^{ \alpha, p } n } { n^{ 1 / \alpha } } \leq x
                \bigg\}
          - G_{\alpha, p, \gamma_n}(x)
            + \frac{1}{2n} G_{\alpha, p, \gamma_n}^{(0,2)}(x)
            - \frac{\mu_1^{\alpha,p}}{n^{1 / \alpha}}
              G_{\alpha, p, \gamma_n}^{(1,1)}(x) \Bigg|
            = o \bigg( \frac{1}{n^{ 1 / \alpha } } \bigg) ,
 \]
 where, for \ $ \alpha \in (0,2) $, \ $ p \in (0,1) $, \ $ \gamma \in (q,1] $
 \ and \ $ k, j \in \{ 0, 1, 2, \dots\} $, \ the function
 \ $G_{\alpha, p, \gamma}^{(k,j)}$ \ can be given as a function of bounded variation
 on the whole \ $ \RR $ \ with Fourier--Stieltjes transform  
 \begin{equation} \label{g^kj}
   \bg_{\alpha, p, \gamma}^{(k,j)}(t)
   = \int_{-\infty}^\infty \ee^{\ii t x} \, \dd G_{\alpha, p, \gamma}^{(k,j)}(x)
   = (- \ii t)^k \, [y_\gamma^{\alpha, p}(t)]^j \, \ee^{ y_\gamma^{\alpha, p}(t) } ,
   \qquad t \in \RR ,
 \end{equation}
 satisfying \ $\lim_{x \to - \infty} G_{\alpha, p, \gamma}^{(k,j)}(x) = 0$.
\ The function \ $G_{\alpha, p, \gamma}^{(k,j)}$ \ can also be given by
 \[
   G_{\alpha, p, \gamma}^{(k,j)}(x)
   = \left. \frac{\partial^{k+j} G_{\alpha, p, \gamma}(x,u)}
                 {\partial x^k \, \partial u^j} \right|_{u=1} , \qquad x \in \RR,
 \]
 where, for each \ $ u > 0 $, \ the function
 \ $ x \mapsto G_{\alpha, p, \gamma}(x,u) $ \ is a semistable infinitely divisible
 distribution function given by its characteristic function
 \ $ t \mapsto \bg_{\alpha, p, \gamma}(t,u) $ \ defined by
 \[
   \bg_{\alpha, p, \gamma}(t,u)
   = \int_{-\infty}^\infty \ee^{\ii t x} \, \dd G_{\alpha, p, \gamma}(x,u)
   = \ee^{ u \, y_\gamma^{\alpha, p}(t) } ,
   \qquad t \in \RR ,
 \]
 see Cs\"org\H{o} \cite{Cso_07b}.

The aim of the present paper is to establish asymptotic expansions of arbitrary
 length for \ $\alpha \in (0,1) \cup (1,2)$, \ improving the approximation.
First, in order to establish candidates for the approximating functions, a
 formal infinite expansion of the characteristic functions of the cumulative
 winnings is derived in Section 2, which can be turned into a formal infinite
 expansion of the distribution functions of the cumulative winnings.
Next, in Section 3, a uniform bound is given in case
 \ $ \alpha \in (0,1) \cup (1,2)$.
\ Its proof, carried out in Section 4, is based on Esseen's classical lemma.
In Section 5 a nonuniform bound is presented for \ $ \alpha \in (1,2) $. 
\ Its proof in Section 6 is based on a lemma of Osipov \cite{Osi_72}.
The convergence rate of the expansion of length \ $\ell$ \ is
 \ $O \left( \frac{1}{n^{1/\alpha}} \right)$ \ for
 \ $\alpha \in \left[ \frac{1}{\ell} , 1 \right)
               \cup \left( 1, 2 - \frac{1}{\ell} \right]$.
\ In the non-lattice case the remainder
 \ $O \left( \frac{1}{n^{1/\alpha}} \right)$ \ will be reduced to
 \ $o \left( \frac{1}{n^{1/\alpha}} \right)$.
\ These convergence rates seem to be optimal. 
The optimality might be studied by the "leading term" approach of Hall
 \cite{Hall_82}.
In case \ $\alpha = 1$ \ one can not improve the short expansion of length
 \ $2$ \ obtained by Cs\"org\H{o} \cite{Cso_07b} by taking a longer expansion,
 since all the extra terms of a formal infinite expansion are of smaller order
 than the remainder term \ $O \left( \frac{1}{n} \right)$. 
\ We also give asymptotic expansions of length \ $\ell$ \ in local merging
 theorems with uniform bound for \ $ \alpha \in (0,1) $ \ and \ $\alpha = 1$,
 \ and with nonuniform bound for \ $ \alpha \in (1,2) $, \ with convergence
 rates \ $O \left( \frac{1}{n^\ell} \right)$,
 \ $O \left( \frac{[\log_r n]^{2 \ell}}{n^\ell} \right)$ \ and
 \ $O \left( \frac{1}{n^{\ell(2-\alpha)/\alpha}} \right)$, \ respectively.

Recently Kevei \cite{Kev_09} derived uniform merging asymptotic expansions for
 distribution functions from the domain of geometric partial attraction of a
 semistable law.

\section{Construction of a formal infinite expansion}
\label{formal_infinite_expansion}

For fixed \ $ \alpha \in (0,2) $ \ and \ $ p \in (0,1) $, \ consider the
 suitably centered and normalized cumulative winnings
 \ $(S_n - c_n^{\alpha, p}) / n^{1 / \alpha}$.
\ We derive a formal infinite expansion (not taking into consideration its
 convergence) for the characteristic function
 \[
   \bff_n^{\alpha, p}(t)
   := \EE \left( \ee^{ \ii t (S_n - c_n^{\alpha, p} ) / n^{1 / \alpha} }
          \right) , \qquad t \in \RR ,
 \]
 in terms of the Fourier--Stieltjes transforms \ $ \bg_{\alpha, p, \gamma_n}^{(k,j)} $,
 \ $ k, j \in \{ 0, 1, 2, \dots \} $.
\ By \eqref{f},
 \begin{equation} \label{f_n}
  \bff_n^{\alpha, p}(t)
  = \big[ \bff_{\alpha,p}( t / n^{1/\alpha} ) \big]^n
    \ee^{ - \ii t c_n^{\alpha, p} / n^{1 / \alpha} }
  = \left[ 1 + \frac{x_n^{\alpha, p}(t)}{n} \right]^n \!
     \ee^{ - \ii t c_n^{\alpha, p} / n^{1 / \alpha} }
 \end{equation}
 for \ $ t \in\RR $, \ where
  \[ 
    x_n^{\alpha, p}(t)
    := n \big( \bff_{\alpha,p}( t / n^{1/\alpha} ) -1 \big)
    = n \sum_{k=1}^\infty
         \Bigg( \exp \bigg\{ \frac{\ii t r^{k / \alpha}}{ n^{1 /\alpha} } \bigg\}
                - 1 \Bigg)
         q^{k-1} p .
  \]
It is easy to check that
 \begin{equation} \label{x}
   x_n^{\alpha, p}(t)
   = y_{\gamma_n}^{\alpha, p}(t) + \ii t \frac{c_n^{\alpha, p}}{ n^{1 / \alpha} }
     + R_{n,1,2}^{\alpha, p}(t) , \qquad t \in \RR ,
 \end{equation}
 where the position parameter \ $\gamma_n$ \ is given in \eqref{gamma_n}, the
 function \ $y_\gamma^{\alpha, p}$ \ stands in the exponent in the definition
 \eqref{g_apg} of the characteristic function \ $\bg_{\alpha, p, \gamma}$, \ and for
 \ $k \in \{2,3,\dots\}$ \ and \ $t \in \RR$,
 \begin{equation} \label{R1k}
   R_{n,1,k}^{\alpha,p}(t)
   := - \sum_{m = \lceil \log_r n \rceil}^\infty
         \left( \exp \left\{ \frac{\ii t}{r^{m/\alpha} \gamma_n^{1/\alpha}} \right\}
                - \sum_{j=0}^{k-1}
                   \frac{(\ii t)^j}
                        {j! \, r^{jm/\alpha} \, \gamma_n^{j/\alpha}} \right)
         \frac{p \gamma_n}{q} r^m ,
 \end{equation}
 see Cs\"org\H{o} \cite[page 833]{Cso_02}.
This series is absolutely convergent on the whole line \ $\RR$ \ for each
 \ $k \geq 2$ \ (see Cs\"org\H{o} \cite[page 322]{Cso_07b} for \ $k = 2$), thus
 \begin{align}
   R_{n,1,k}^{\alpha, p}(t)
   & = - \frac{p \gamma_n}{q}
         \sum_{m = \lceil \log_r n \rceil}^\infty \,
          \sum_{j=k}^\infty \,
           \frac{(\ii t)^j \, r^m}{j! \, r^{jm/\alpha} \, \gamma_n^{j/\alpha}}
     \nonumber \\[2mm]
   & = - \frac{p}{q}
         \sum_{j=k}^\infty
          \frac{(\ii t)^j}{j! \, \gamma_n^{(j - \alpha) / \alpha}}
          \sum_{m = \lceil \log_r n \rceil}^\infty q^{m(j - \alpha) / \alpha} 
     \nonumber \\[2mm]
   & = - \frac{p}{q}
         \sum_{j=k}^\infty
          \frac{(\ii t)^j \, q^{\lceil \log_r n \rceil (j - \alpha) / \alpha}}
               {j! \, \gamma_n^{(j - \alpha) / \alpha} \, ( 1 - q^{(j - \alpha) / \alpha} )}
     \nonumber \\[2mm]
   & = - \frac{p}{q}
         \sum_{j=k}^\infty
          \frac{(\ii t)^j}
               {j! \, ( 1 - q^{(j - \alpha) / \alpha} ) \, n^{(j - \alpha) / \alpha} }
     = n \sum_{j=k}^\infty
          \frac{\mu_j^{\alpha,p}}{j!}
          \left( \frac{\ii t}{n^{1 / \alpha}} \right)^j . \label{R1}
 \end{align}
Consequently, we have the formal expansion 
 \begin{align*}  
  \bff_n^{\alpha, p}(t)
  &= \exp \Bigg\{ - \ii t \frac{c_n^{\alpha, p}}{ n^{1 / \alpha} }
                  + n \log \bigg( 1 + \frac{x_n^{\alpha,p}(t)}{n} \bigg)
          \Bigg\} 
     \\[2mm]
  &= \exp \Bigg\{ - \ii t \frac{c_n^{\alpha, p}}{ n^{1 / \alpha} }
                  + x_n^{\alpha,p}(t)
                  + n \sum_{k=2}^\infty 
                       \frac{(-1)^{k+1}}{k}
                       \left( \frac{x_n^{\alpha,p}(t)}{n} \right)^k \Bigg\} 
     \\[2mm]  
  &= \exp \Bigg\{ y_{\gamma_n}^{\alpha,p}(t)
                  + n \sum_{j=2}^\infty
                       \frac{\mu_j^{\alpha,p}}{j!}
                       \left( \frac{\ii t}{n^{1 / \alpha}} \right)^j
                  + n \sum_{k=2}^\infty 
                       \frac{(-1)^{k+1}}{k}
                       \left( \frac{x_n^{\alpha,p}(t)}{n} \right)^k \Bigg\} . 
 \end{align*}
First consider the case \ $ \alpha \in (0,1) \cup (1,2) $. 
\ Then by \eqref{x} and \eqref{c}, we have
 \begin{align*}   
  \bff_n^{\alpha, p}(t)
  &= \exp \Bigg\{ y_{\gamma_n}^{\alpha,p}(t)
                  + n \sum_{j=2}^\infty
                       \frac{\mu_j^{\alpha,p}}{j!}
                       \left( \frac{\ii t}{n^{1 / \alpha}} \right)^j +
     \\[2mm]
  &\phantom{\quad\;\exp \Bigg\{}
                  + n \sum_{k=2}^\infty 
                       \frac{(-1)^{k+1}}{k}
                       \left[ \frac{y_{\gamma_n}^{\alpha,p}(t)}{n}
                              + \sum_{\ell=1}^\infty
                                 \frac{\mu_\ell^{\alpha,p}}{\ell!}
                                 \left( \frac{\ii t}{n^{1 / \alpha}} \right)^\ell
                       \right]^k \Bigg\} 
     \\[2mm]
  &= \ee^{y_{\gamma_n}^{\alpha,p}(t)}
     \exp \Bigg\{ n \sum_{k=0}^\infty
                     \sum_{j=0}^\infty 
                      d_{k,j}^{\alpha,p}
                       \left( - \frac{\ii t}{n^{1 / \alpha}} \right)^k
                       \bigg( \frac{y_{\gamma_n}^{\alpha,p}(t)}{n} \bigg)^j
          \Bigg\} , 
 \end{align*}
 where the coefficients \ $d_{k,j}^{\alpha,p}$ \ are polynomials of the virtual
 moments \ $\mu_\ell^{\alpha,p}$, \ $\ell \in \NN$, \ such that
 \ $d_{0,0}^{\alpha,p} = d_{0,1}^{\alpha,p} = d_{1,0}^{\alpha,p} = 0$.
\ Remark that
 \begin{gather*}
  d_{0,2}^{\alpha,p} = - \frac{1}{2} , \qquad
   d_{1,1}^{\alpha,p} = \mu_1^{\alpha,p} , \qquad
   d_{2,0}^{\alpha,p} = - \frac{1}{2} [\mu_1^{\alpha,p}]^2
                    + \frac{1}{2} \mu_2^{\alpha,p} , \\[2mm]
  d_{0,3}^{\alpha,p} = \frac{1}{3} , \qquad
   d_{1,2}^{\alpha,p} = - \mu_1^{\alpha,p} , \qquad
   d_{2,1}^{\alpha,p} = [\mu_1^{\alpha,p}]^2 - \frac{1}{2} \mu_2^{\alpha,p} ,
  d_{3,0}^{\alpha,p} = - \frac{1}{3} [\mu_1^{\alpha,p}]^3
                   + \frac{1}{2} \mu_1^{\alpha,p} \mu_2^{\alpha,p} 
                   - \frac{1}{6} \mu_3^{\alpha,p} , \\[2mm]
  d_{0,4}^{\alpha,p} = - \frac{1}{4} , \qquad
   d_{1,3}^{\alpha,p} = \mu_1^{\alpha,p} , \qquad
   d_{2,2}^{\alpha,p} = - \frac{3}{2} [\mu_1^{\alpha,p}]^2
                    + \frac{1}{2} \mu_2^{\alpha,p} , \qquad 
  d_{3,1}^{\alpha,p} = [\mu_1^{\alpha,p}]^3
                   - \mu_1^{\alpha,p} \mu_2^{\alpha,p} 
                   + \frac{1}{6} \mu_3^{\alpha,p} , \\[2mm]
  d_{4,0}^{\alpha,p} = - \frac{1}{4} [\mu_1^{\alpha,p}]^4
                   + \frac{1}{2} [\mu_1^{\alpha,p}]^2 \mu_2^{\alpha,p} 
                   - \frac{1}{6} \mu_1^{\alpha,p} \mu_3^{\alpha,p}
                   - \frac{1}{8} [\mu_2^{\alpha,p}]^2
                   + \frac{1}{24} \mu_4^{\alpha,p} .
 \end{gather*}
Following Christoph and Wolf \cite[Section 4.3]{Cri_Wolf_92}, using the Taylor
 expansion of the exponential function, we obtain a formal power series
 expansion
 \[
   \bff_n^{\alpha, p}(t)
   = \ee^{y_{\gamma_n}^{\alpha,p}(t)}
     \Bigg\{ 1 + \sum_{k=0}^\infty \,
                  \sum_{j = - \lfloor k/2 \rfloor}^\infty \,
                   \sum_{m = \max\{ 1, -j \}}^{k+j}
                    \frac{w_{m,k,j}^{\alpha, p} \, (-\ii t)^k \,
                          [y_{\gamma_n}^{\alpha,p}(t)]^{j+m}}
                         {m! \, n^{\frac{k}{\alpha} + j}} \Bigg\},
 \] 
 where
 \[
   w_{m,k,j}^{\alpha, p}
   := \begin{cases}
       0, & \text{if \ $- \lfloor k/2 \rfloor \leq j < L_{m,k}$,} \\[1mm]
       \sideset{}{'}\sum\limits_{\overset{\hbox{$\SC k_1 + \cdots + k_m = k$}}
                            {\hbox{$\SC s_1 + \cdots + s_m = j + m$}}}
        d_{k_1,s_1}^{\alpha, p} \dots d_{k_m,s_m}^{\alpha, p} ,
        & \text{if \ $j \geq L_{m,k}$,}  
      \end{cases}
 \]
 with
 \[
   L_{m,k} := \max\{ - \lfloor k/2 \rfloor , - m, m - k \}.
 \]
The summation \ $\sideset{}{'}\sum$ \ in the definition of
 \ $ w_{m,k,j}^{\alpha, p}$ \ is carried over all nonnegative integers \ $k_1$,
 \dots,  $k_m$, $s_1$, \dots,  $s_m$, \ such that \ $k_1 + \cdots + k_m = k$,
 \ $s_1 + \cdots + s_m = j + m$, \ and
 \ $(k_j,s_j) \notin \{ (0,0), (1,0), (0,1) \}$.
\ Finally, we conclude a formal expansion
 \begin{equation} \label{infinite_expantion}
   \bff_n^{\alpha, p}(t)
   = \bg_{\alpha, p, \gamma_n}(t)
     + \sum_{k=0}^\infty \, 
        \sum_{j = - \lfloor k/2 \rfloor}^\infty \,
         \sum_{m = \max\{ 1, -j \}}^{k+j}
          \frac{w_{m,k,j}^{\alpha, p} \, \bg_{\alpha, p, \gamma_n}^{(k,\,j+m)}(t)}
               {m! \, n^{\frac{k}{\alpha} + j}} ,
 \end{equation}
 where the position parameter \ $\gamma_n$ \ is introduced in \eqref{gamma_n},
 the characteristic function \ $\bg_{\alpha, p, \gamma}$ \ is defined in
 \eqref{g_apg}, and the function \ $\bg_{\alpha, p, \gamma}^{(k,\,j)}$ \ is given in
 \eqref{g^kj}.
Replacing each term by its inverse Fourier-Stieltjes transform, this may be
 turned into a formal expansion
 \[
   \PP \! \! \bigg\{ \! \frac{S_n \! - \! \mu_1^{ \alpha, p } n }
                             { n^{ 1 / \alpha } } \! \leq x \! \bigg\} \!
   = G_{\alpha, p, \gamma_n}(x)
     + \sum_{k=0}^\infty \,
        \sum_{j = - \lfloor k/2 \rfloor}^\infty \,
         \sum_{m = \max\{ 1, -j \}}^{k+j} \hspace*{-3mm}
          \frac{w_{m,k,j}^{\alpha, p} \, G_{\alpha, p, \gamma_n}^{(k,\,j+m)}(x)}
               {m! \, n^{\frac{k}{\alpha} + j}} ,
 \] 
 for \ $x \in\RR$.

Now consider the case \ $ \alpha = 1$. 
\ Then by \eqref{x} and \eqref{c}, we have
 \begin{align*}   
  \bff_n^{1, p}(t)
  &= \exp \Bigg\{ y_{\gamma_n}^{1,p}(t)
                  + n \sum_{j=2}^\infty
                       \frac{\mu_j^{1,p}}{j!}
                       \left( \frac{\ii t}{n} \right)^j
     \\[2mm]
  &\phantom{\quad\;\exp \Bigg\{}
                  + n \sum_{k=2}^\infty 
                       \frac{(-1)^{k+1}}{k}
                       \left[ \frac{y_{\gamma_n}^{1,p}(t) + \ii t p r \log_r n}{n}
                              + \sum_{\ell=2}^\infty
                                 \frac{\mu_\ell^{1,p}}{\ell!}
                                 \left( \frac{\ii t}{n} \right)^\ell
                       \right]^k \Bigg\} 
     \\[2mm]
  &= \ee^{y_{\gamma_n}^{1,p}(t)}
     \exp \Bigg\{ n \sum_{k=0}^\infty
                     \sum_{j=0}^\infty 
                      d_{k,j}^{1,p}
                       \left( - \frac{\ii t}{n} \right)^k
                       \bigg( \frac{y_{\gamma_n}^{1,p}(t) + \ii t p r \log_r n}{n}
                       \bigg)^j
          \Bigg\} , 
 \end{align*}
 where the coefficients \ $d_{k,j}^{1,p}$ \ are now polynomials of the virtual
 moments \ $\mu_\ell^{1,p}$ \ with only \ $\ell \geq 2$, \ such that
 \ $d_{0,0}^{1,p} = d_{0,1}^{1,p} = d_{1,0}^{1,p} = 0$.
\ Remark that a formula for \ $d_{k,j}^{1,p}$ \ can be obtained from the formula
 for \ $d_{k,j}^{\alpha,p}$, \ $\alpha \in (0,1) \cup (1,2)$, \ by replacing
 \ $\alpha$ \ by \ $1$ \ and \ $\mu_1^{1,p}$ \ by \ $0$.
\ Hence following Christoph and Wolf \cite[Section 4.3]{Cri_Wolf_92}, now we
 obtain a formal power series expansion
 \[
   \ee^{y_{\gamma_n}^{1,p}(t)}
   \Bigg\{ 1 + \sum_{k=0}^\infty \,
                \sum_{j = - \lfloor k/2 \rfloor}^\infty \,
                 \sum_{m = \max\{ 1, -j \}}^{k+j}
                  \frac{w_{m,k,j}^{1, p} \, (-\ii t)^k \,
                        [y_{\gamma_n}^{1,p}(t) + \ii t p r \log_r n]^{j+m}}
                       {m! \, n^{k + j}} \Bigg\} 
 \]
 for \ $\bff_n^{1, p}(t)$, \ $t \in \RR$, \ which may be turned into a formal
 expansion
 \begin{multline*}
  \PP \bigg\{ \frac{S_n - p r \, n \log_r n}{ n } \leq x \bigg\}
  = G_{1, p, \gamma_n}(x)
    + \sum_{k=0}^\infty \,
        \sum_{j = - \lfloor k/2 \rfloor}^\infty \,
         \sum_{m = \max\{ 1, -j \}}^{k+j}
          \frac{w_{m,k,j}^{1, p}}{m! \, n^{k + j}} \\[2mm]
          \times
          \sum_{\ell = 0}^{j+m}
           \binom{j+m}{\ell}
           [- p r \log_r n]^\ell \,
           G_{1, p, \gamma_n}^{(k+\ell,\,j+m-\ell)}(x) ,
 \end{multline*}
 for \ $x \in \RR$.
\ On the other hand, the above formula for \ $\bff_n^{1, p}(t)$, \ and hence the
 above expansion for
 \ $\PP \left\{ \frac{S_n - p r \, n \log_r n}{ n } \leq x \right\}$ \ in case 
 \ $\alpha = 1$ \ can be obtained from the formula
 for \ $\bff_n^{\alpha, p}(t)$, \ $\alpha \in (0,1) \cup (1,2)$, \ and from the
 expansion for 
 \ $\PP \left\{ \frac{S_n - \mu_1^{ \alpha, p } n } { n^{ 1 / \alpha } } \leq x
        \right\}$
 \ by replacing \ $\alpha$ \ by \ $1$ \ and \ $\mu_1^{1,p}$ \ by
 \ $p r \log_r n$.

\section{Uniform bounds in asymptotic expansions}
\label{unif_bounds}

For \ $\alpha \in (0,1) \cup (1,2)$, \ $p \in (0,1)$, \ $n \in \NN$ \ and
 \ $\ell \in \{ 0, 1, \dots\}$, \ introduce
 \[
   G_{n,\ell}^{\alpha,p}(x)
   := G_{\alpha, p, \gamma_n}(x)
      + \sum_{k=0}^{2 \ell} \,
        \sum_{j = - \lfloor k/2 \rfloor}^{\ell-k} \,
         \sum_{m = \max\{ 1, - j \}}^{k+j} 
          \frac{w_{m,k,j}^{\alpha, p} \, G_{\alpha, p, \gamma_n}^{(k,j+m)}(x)}
               {m!  \, n^{\frac{k}{\alpha} + j}}
 \]
 for \ $x \in \RR$.
\ The function \ $G_{n,\ell}^{\alpha,p}$ \ consists of all terms of order
 \ $n^{- \frac{k}{\alpha} - j}$ \ with \ $k + j \leq \ell$ \ of the formal expansion. 
Remark that \ $G_{n,0}^{\alpha,p} = G_{\alpha, p, \gamma_n}$, \ and
 \begin{align*}
  G_{n,1}^{\alpha,p}
  & = G_{\alpha, p, \gamma_n}
      + \frac{d_{0,2}^{\alpha,p} \, G_{\alpha, p, \gamma_n}^{(0,2)}}{n}  
      + \frac{d_{1,1}^{\alpha,p} \, G_{\alpha, p, \gamma_n}^{(1,1)}}{n^{\frac{1}{\alpha}}}
      + \frac{d_{2,0}^{\alpha,p} \, G_{\alpha, p, \gamma_n}^{(2,0)}}
             {n^{\frac{2}{\alpha} - 1}} \\[2mm]
  & = G_{\alpha, p, \gamma_n}
      - \frac{G_{\alpha, p, \gamma_n}^{(0,2)}}{2 n}  
      + \frac{\mu_1^{\alpha,p} \, G_{\alpha, p, \gamma_n}^{(1,1)}}{n^{\frac{1}{\alpha}}}
      - \frac{\big( [\mu_1^{\alpha,p}]^2 - \mu_2^{\alpha,p} \big) \,
              G_{\alpha, p, \gamma_n}^{(2,0)}}
             {2 n^{\frac{2}{\alpha} - 1}} ,
 \end{align*}
 \begin{align*}
  & G_{n,2}^{\alpha,p}
    = G_{n,1}^{\alpha,p}
      + \frac{d_{0,3}^{\alpha,p} \, G_{\alpha, p, \gamma_n}^{(0,3)} 
              + \frac{1}{2} [d_{0,2}^{\alpha,p}]^2 \, G_{\alpha, p, \gamma_n}^{(0,4)}}
             {n^2}
      + \frac{d_{1,2}^{\alpha,p} \, G_{\alpha, p, \gamma_n}^{(1,2)}
              + d_{0,2}^{\alpha,p} \, d_{1,1}^{\alpha,p} \, G_{\alpha, p, \gamma_n}^{(1,3)}}
             {n^{1 + \frac{1}{\alpha}}} \\[2mm] 
  &\phantom{G_{n,2}^{\alpha,p}\quad}
      + \frac{d_{2,1}^{\alpha,p} \, G_{\alpha, p, \gamma_n}^{(2,1)}
              + \frac{1}{2}
                \big( [d_{1,1}^{\alpha,p}]^2
                       + 2 d_{2,0}^{\alpha,p} d_{0,2}^{\alpha,p} \big) \,
                G_{\alpha, p, \gamma_n}^{(2,2)}}
             {n^{\frac{2}{\alpha}}} \\[2mm] 
  &\phantom{G_{n,2}^{\alpha,p}\quad}
      + \frac{d_{3,0}^{\alpha,p} \, G_{\alpha, p, \gamma_n}^{(3,0)}
              + d_{1,1}^{\alpha,p} \, d_{2,0}^{\alpha,p} \, G_{\alpha, p, \gamma_n}^{(3,1)}}
             {n^{\frac{3}{\alpha} - 1}}
      + \frac{\frac{1}{2} [d_{2,0}^{\alpha,p}]^2 \, G_{\alpha, p, \gamma_n}^{(4,0)}}
             {n^{\frac{4}{\alpha} - 2}} \\[2mm] 
  &\phantom{G_{n,2}^{\alpha,p}}
     = G_{n,1}^{\alpha,p}
      + \frac{8 G_{\alpha, p, \gamma_n}^{(0,3)} + 3 G_{\alpha, p, \gamma_n}^{(0,4)}}
             {24 \, n^2}
      - \frac{\mu_1^{\alpha,p}
              \big(2 G_{\alpha, p, \gamma_n}^{(1,2)}
                   + G_{\alpha, p, \gamma_n}^{(1,3)}\big)}
             {2 n^{1 + \frac{1}{\alpha}}} \\[2mm] 
  &\phantom{G_{n,2}^{\alpha,p}\quad}
      + \frac{2 \big( 2 [\mu_1^{\alpha,p}]^2 - \mu_2^{\alpha,p} \big) \,
              G_{\alpha, p, \gamma_n}^{(2,1)}
              + \big( 3 [\mu_1^{\alpha,p}]^2 - \mu_2^{\alpha,p} \big) \,
                G_{\alpha, p, \gamma_n}^{(2,2)} }
             {4 n^{\frac{2}{\alpha}}} \\[2mm] 
  &\phantom{G_{n,2}^{\alpha,p}\quad}
      - \frac{ \big( 2 [\mu_1^{\alpha,p}]^3 - 3 \mu_1^{\alpha,p} \, \mu_2^{\alpha,p}
                     + \mu_3^{\alpha,p} \big) \,
               G_{\alpha, p, \gamma_n}^{(3,0)}
              + 3 \mu_1^{\alpha,p}
                  \big( [\mu_1^{\alpha,p}]^2 - \mu_2^{\alpha,p} \big) \,
                  G_{\alpha, p, \gamma_n}^{(3,1)}}
             {6 n^{\frac{3}{\alpha} - 1}} \\[2mm] 
  &\phantom{G_{n,2}^{\alpha,p}\quad}
      + \frac{ \big( [\mu_1^{\alpha,p}]^2 - \mu_2^{\alpha,p} \big)^2 \,
               G_{\alpha, p, \gamma_n}^{(4,0)}}
             {8 n^{\frac{4}{\alpha} - 2}} .
 \end{align*}
For \ $\alpha = 1$, \ $p \in (0,1)$, \ $n \in \NN$ \ and
 \ $\ell \in \{ 0, 1, \dots\}$, \ introduce
 \begin{multline*}
  G_{n,\ell}^{1, p}(x)
  := G_{1, p, \gamma_n}(x)
     + \sum_{k=0}^{2 \ell} \,
        \sum_{j = - \lfloor k/2 \rfloor}^{\ell - k} \,
         \sum_{m = \max\{ 1, -j \}}^{k+j}
          \frac{w_{m,k,j}^{1, p}}{m! \, n^{k + j}} \\[2mm]
          \times
          \sum_{\ell = 0}^{j+m}
           \binom{j+m}{\ell}
           [- p r \log_r n]^\ell \,
           G_{1, p, \gamma_n}^{(k+\ell,\,j+m-\ell)}(x) ,
 \end{multline*}
 for \ $x \in \RR$.
\ The function \ $G_{n,\ell}^{1,p}$ \ consists of all terms containing
 \ $n^{- k - j}$ \ with \ $k + j \leq \ell$ \ of the formal expansion. 
Remark that \ $G_{n,0}^{1,p} = G_{1, p, \gamma_n}$, \ and
 \[
  G_{n,1}^{1,p}
  = G_{1, p, \gamma_n}
      - \frac{G_{1, p, \gamma_n}^{(0,2)}}{2 n}  
      + \frac{[p r \log_r n] \, G_{1, p, \gamma_n}^{(1,1)}}{n}
      - \frac{\big( [p r \log_r n]^2 - \mu_2^{1,p} \big) \, G_{1, p, \gamma_n}^{(2,0)}}
             {2 n} ,
 \]
 \begin{align*}
  G_{n,2}^{1,p}
  & = G_{n,1}^{1,p}
      + \frac{8 G_{1, p, \gamma_n}^{(0,3)} + 3 G_{1, p, \gamma_n}^{(0,4)}}{24 \, n^2}
      - \frac{[p r \log_r n]
              \big(2 G_{1, p, \gamma_n}^{(1,2)} + G_{1, p, \gamma_n}^{(1,3)}\big)}
             {2 n^2} \\[2mm] 
  &\phantom{\quad}
      + \frac{2 \big( 2 [p r \log_r n]^2 - \mu_2^{1,p} \big) \,
              G_{1, p, \gamma_n}^{(2,1)}
              + \big( 3 [p r \log_r n]^2 - \mu_2^{1,p} \big) \,
                G_{1, p, \gamma_n}^{(2,2)} }
             {4 n^2} \\[2mm] 
  &\phantom{\quad}
      - \frac{ \big( 2 [p r \log_r n]^3 - 3 [p r \log_r n] \, \mu_2^{1,p}
                     + \mu_3^{1,p} \big) \,
               G_{1, p, \gamma_n}^{(3,0)}}
             {6 n^2} \\[2mm]  
  &\phantom{\quad}
      - \frac{ 3 [p r \log_r n]
                  \big( [p r \log_r n]^2 - \mu_2^{1,p} \big) \,
                  G_{1, p, \gamma_n}^{(3,1)}}
             {6 n^2}
      + \frac{ \big( [p r \log_r n]^2 - \mu_2^{1,p} \big)^2 \,
               G_{1, p, \gamma_n}^{(4,0)}}
             {8 n^2} .
 \end{align*}
The functions \ $G_{n,1}^{\alpha, p}$, \ $\alpha \in (0,2)$, \ are almost the same
 as the approximating functions \ $G_n^{\alpha, p}$, \ $\alpha \in (0,2)$, \ of
 Cs\"org\H{o} \cite{Cso_07b}. 
They differ only in the coefficients of their last terms, which does not make
 difference in the convergence rates if \ $\alpha \in(0,1]$, \ but the
 functions \ $G_{n,1}^{\alpha, p}$, \ $\alpha \in (1,2)$, \ give better
 approximations than \ $G_n^{\alpha, p}$, \ $\alpha \in (1,2)$, \ of Cs\"org\H{o}
 \cite{Cso_07b}.

The main results of this section are contained in the following

\begin{pro} \label{main_uniform}
For \ $\ell \in \NN$,
 \[
   \sup_{ x \in \RR }
    \Bigg| \! \PP \! \! \bigg\{ \! \frac{S_n \! - \mu_1^{ \alpha, p } n }
                                        { n^{ 1 / \alpha } } \!
                                \leq x \! \bigg\}
           - G_{n,\ell - 1}^{\alpha, p}(x) \Bigg|
   = \begin{cases}
      O \! \left( \frac{1}{n^\ell} \right) ,
       & \text{if \ $ 0 < \alpha < \frac{1}{\ell} $,} \\[2mm]
      O \! \left( \frac{1}{n^{1 / \alpha}} \right) ,
       & \text{if \ $ \frac{1}{\ell} \leq \alpha < 1 $
               \ or \ $ 1 < \alpha \leq 2 - \frac{1}{\ell} $,} \\[2mm]
      O \! \left( \frac{1}{ n^{ \ell ( 2 - \alpha ) / \alpha } } \right) ,
       & \text{if \ $ 2 - \frac{1}{\ell} < \alpha < 2 $.}
     \end{cases}
 \] 
For \ $\ell \in \{ 2, 3, \dots\}$,
 \ $ \alpha \in \left( \frac{1}{\ell} , 1 \right)
                \cup \left( 1, 2 - \frac{1}{\ell} \right) $
 \ and \ $ r^{ 1 / \alpha} \notin \NN $,
 \[
   \sup_{ x \in \RR }
    \Bigg| \PP \! \bigg\{ \frac{S_n \! - \mu_1^{ \alpha, p } n }
                               { n^{ 1 / \alpha } } \!
                          \leq x \bigg\}
           - G_{n,\ell - 1}^{\alpha, p}(x) \Bigg|
   = o \left( \frac{1}{n^{1 / \alpha}} \right).
 \] 
\end{pro}

In case \ $\ell = 1$, \ the approximating function is
 \ $G_{n,0}^{\alpha, p} = G_{\alpha, p, \gamma_n}$, \ and Proposition \ref{main_uniform}
 gives back the rates of merge for \ $0 < \alpha < 1$ \ and \ $1 < \alpha < 2$,
 \ due to Cs\"org\H{o} \cite{Cso_02}.
Certain terms of \ $G_{n,\ell - 1}^{\alpha, p}$ \ are of the same or of a smaller
 order than the remainder terms \ $O \left( \frac{1}{n^{1/\alpha}} \right)$ \ or
 \ $o \left( \frac{1}{n^{1/\alpha}} \right)$.
\ Using boundedness of the functions \ $G_{\alpha, p, \gamma_n}^{(k,j)}$ \ (see 
 Cs\"org\H{o} \cite[Lemma 6]{Cso_07b}), the expansions may be simplified as
 follows.

\begin{thm} \label{MAIN_uniform} 
For \ $\ell \in \{ 2, 3, \dots\}$,
 \ $ \alpha \in \left[ \frac{1}{\ell} , \frac{1}{\ell - 1} \right)
                \cup \left( 2 - \frac{1}{\ell-1} , \,
                            2 - \frac{1}{\ell} \right] $
 \ and \ $ r^{ 1 / \alpha} \in \NN $,
 \[
   \sup_{ x \in \RR }
    \Bigg| \PP \bigg\{ \frac{S_n \! - \mu_1^{ \alpha, p } n }
                            { n^{ 1 / \alpha } } \!
                       \leq x \bigg\}
           - \tG_{n,\ell - 1}^{\alpha, p}(x) \Bigg|
   = O \left( \frac{1}{n^{1 / \alpha}} \right) ;
 \]
 for \ $\ell = 2$,
 \ $\alpha \in \left( \frac{1}{2} , 1 \right)
                      \cup \left( 1, \frac{3}{2} \right)$,
 \ or for \ $\ell \in \{ 3, 4, \dots\}$,
 \ $ \alpha \in \left( \frac{1}{\ell} , \frac{1}{\ell - 1} \right]
                \cup \left[ 2 - \frac{1}{\ell-1} , \,
                            2 - \frac{1}{\ell} \right)$,
 \ and for \ $ r^{ 1 / \alpha} \notin \NN $,
 \[
   \sup_{ x \in \RR }
    \Bigg| \PP \bigg\{ \frac{S_n \! - \mu_1^{ \alpha, p } n }
                            { n^{ 1 / \alpha } } \!
                       \leq x \bigg\}
           - \ttG_{n,\ell - 1}^{\alpha, p}(x) \Bigg|
   = o \left( \frac{1}{n^{1 / \alpha}} \right) ,
 \]
 where the approximating functions are given by
 \[
   \tG_{n,\ell - 1}^{\alpha, p}
   := \begin{cases}
       \DS G_{\alpha, p, \gamma_n}
       + \sum_{j=1}^{\ell - 1}
          \sum_{m=1}^j
           \frac{w_{m,0,j}^{\alpha, p} \, G_{\alpha, p, \gamma_n}^{(0,j+m)}}
                { m!  \, n^j } ,
        & \text{if \ $ \alpha \in (0, 1) $,} \\[5mm]
       \DS G_{\alpha, p, \gamma_n}
       + \sum_{j=1}^{\ell - 1}
          \frac{(-1)^j \big( [\mu_1^{\alpha,p}]^2 - \mu_2^{\alpha,p} \big)^j \,
                G_{\alpha, p, \gamma_n}^{(2j,0)}}
               { m! \, 2^j \, n^{ j ( 2 - \alpha ) / \alpha } } ,
        & \text{if \ $ \alpha \in (1, 2) $,}
      \end{cases}
 \]
 and
 \ $\DS \ttG_{n,\ell - 1}^{\alpha, p}
    := \tG_{n,\ell - 1}^{\alpha, p}
       + \frac{\mu_1^{\alpha,p} \, G_{\alpha, p, \gamma_n}^{(1,1)}}{n^{1 / \alpha}}$.
\end{thm}

If \ $\ell = 2$ \ then the approximating function is
 \[
   \tG_{n,1}^{\alpha, p}
   = \begin{cases}
      \DS G_{\alpha, p, \gamma_n} - \frac{G_{\alpha, p, \gamma_n}^{(0,2)}}{2n} ,
       & \text{if \ $ \alpha \in (0, 1) $,} \\[2mm]
      \DS G_{\alpha, p, \gamma_n}
      - \frac{\big( [\mu_1^{\alpha,p}]^2 - \mu_2^{\alpha,p} \big) \,
              G_{\alpha, p, \gamma_n}^{(2,0)}}
             { 2 \, n^{ ( 2 - \alpha ) / \alpha } } ,
        & \text{if \ $ \alpha \in (1, 2) $,}
     \end{cases}
 \]
 and for \ $\alpha \in (0, 1)$ \ the result has been obtained by Cs\"org\H{o}
 \cite{Cso_07b}.
The approximating functions for \ $\ell = 3$, \ $\ell = 4$ \ and
 \ $\alpha \in (0, 1)$ \ are
 \begin{align*}
   \tG_{n,2}^{\alpha, p}
   & = \tG_{n,1}^{\alpha, p}
       + \frac{8 G_{\alpha, p, \gamma_n}^{(0,3)} + 3 G_{\alpha, p, \gamma_n}^{(0,4)}}
              {24n^2} , \\[2mm]
   \tG_{n,3}^{\alpha, p}
   & = \tG_{n,2}^{\alpha, p}
       - \frac{12 G_{\alpha, p, \gamma_n}^{(0,4)} + 8 G_{\alpha, p, \gamma_n}^{(0,5)}
               + G_{\alpha, p, \gamma_n}^{(0,6)}}
              {48 n^3} .
 \end{align*}

\section{Proof of Proposition \ref{main_uniform}}

Fix \ $\ell \in \NN$ \ and \ $\alpha \in (0,1) \cup (1,2)$.
\ The strategy of the proof is the same as the proof of the Proposition of
 Cs\"org\H{o} \cite{Cso_07b}.
It is based on Esseen's classical result (see Petrov
 \cite[Theorem 5.2]{Pet_75}).

\begin{lem}[Esseen] \label{Esseen1}
Let \ $F$ \ be a distribution function and \ $G$ \ be a differentiable function
 of bounded variation on \ $\RR$ \ with Fourier--Stieltjes transforms
 \ $ \bff(t) = \int_{-\infty}^\infty \ee^{\ii t x} \, \dd F(x)$ \ and
 \ $ \bg(t) = \int_{-\infty}^\infty \ee^{\ii t x} \, \dd G(x)$, \ $t \in\RR$, \ such
 that \ $G(-\infty): = \lim_{x\to-\infty} G(x) = 0$.
\ Then
 \[
   \sup_{x \in \RR} | F(x) - G(x) |
   \leq \frac{b}{2 \pi}
        \int_{-T}^T
         \left| \frac{\bff(t) - \bg(t)}{t} \right| \, \dd t
        + c_b \frac{\sup_{x \in \RR} | G'(x) |}{T}         
 \]
 for every choice of \ $T > 0$ \ and \ $b > 1$, \ where \ $c_b > 0$ \ is a
 constant depending only on \ $b$, \ which can be given as
 \ $c_b = 4 b d_b^2 / \pi$, \ where \ $d_b > 0$ \ is the unique root \ $d$ \ of
 the equation
 \ $\frac{4}{\pi} \int_0^d \frac{\sin^2 u}{u^2} \, \dd u = 1 + \frac{1}{b}$.
\end{lem}

We apply this lemma for
 \ $F(x) = \PP \left\{ \frac{S_n \! - \mu_1^{ \alpha, p } n }
                            { n^{ 1 / \alpha } } \!
                       \leq x \right\}$
 \ and \ $G(x) = G_{n,\ell - 1}^{\alpha, p}(x)$, \ $x \in \RR$.
\ By Lemma 6 of Cs\"org\H{o} \cite{Cso_07b}, \ $G_{n,\ell - 1}^{\alpha, p}$ \ is a
 differentiable function of bounded variation on \ $\RR$ \ with
 \ $G_{n,\ell - 1}^{\alpha, p}(-\infty) = 0$.
\ The Fourier--Stieltjes transform of \ $G_{n,\ell - 1}^{\alpha, p}$ \ is
 \[
   \bg_{n,\ell-1}^{\alpha,p}(t)
   = \ee^{y_{\gamma_n}^{\alpha,p}(t)}
     \Bigg\{ 1 + \sum_{k=0}^{2 \ell-2} \,
                  \sum_{j = - \lfloor k/2 \rfloor}^{\ell-k-1} \,
                   \sum_{m = \max\{ 1, - j \}}^{k+j} \hspace*{-2mm}
                    \frac{w_{m,k,j}^{\alpha, p} \, (- \ii t)^k \,
                          [y_{\gamma_n}^{\alpha,p}(t)]^{j+m}}
                         {m!  \, n^{\frac{k}{\alpha} + j}} \Bigg\} ,
 \]
 for \ $t \in \RR$.
\ The aim of the following discussion is to find an appropriate estimate for
 \ $| \bff_n^{\alpha,p}(t) - \bg_{n,\ell-1}^{\alpha,p}(t) |$ \ if
 \ $|t| \leq T_{n,\ell-1}^{\alpha,p}$ \ with appropriate \ $T_{n,\ell-1}^{\alpha,p}$.
\ Recalling formula \eqref{f_n} for \ $\bff_n^{\alpha,p}$, \ we have
 \begin{equation} \label{f_n^ap}
   \bff_n^{\alpha,p}(t)
   = \exp \left\{ - \ii t \mu_1^{\alpha, p} n^{ ( \alpha - 1 ) / \alpha  }
                  + n \log \left( 1 + \frac{x_n^{\alpha, p}(t)}{n} \right)
          \right\} ,
 \end{equation}
 provided \ $|x_n^{\alpha, p}(t)| < n$.
\ In order to estimate \ $|x_n^{\alpha, p}(t)|$, \ we will use formula \eqref{x}
 for \ $x_n^{\alpha, p}(t)$ \ containing \ $y_{\gamma_n}^{\alpha,p}(t)$ \ and
 \ $R_{n,1,2}^{\alpha,p}(t)$. 
\ Cs\"org\H{o} \cite[Lemma 3]{Cso_02}, \cite[Lemma 3]{Cso_07b} proved the
 following crucial estimates.

\begin{lem}[Cs\"org\H{o}] \label{y}
For arbitrary \ $\alpha \in (0,1) \cup (1,2)$ \ and \ $p \in (0,1)$, \ there
 exists \ $C_1^{\alpha,p} > 0$ \ such that, uniformly in \ $\gamma \in (q,1]$,
 \[
   \Re ( y_\gamma^{\alpha,p}(t) ) \leq - C_1^{\alpha,p} |t|^\alpha, \qquad
   | y_\gamma^{\alpha,p}(t) | \leq C_1^{\alpha,p} |t|^\alpha, \qquad  t \in \RR .
 \]
\end{lem}

Applying this lemma, we can derive the following estimates. 

\begin{lem} \label{R1_x}
For arbitrary \ $k \geq 2$, \ $\alpha \in (0,1) \cup (1,2)$ \ and
 \ $p \in (0,1)$, \ there exists \ $C_{1,k}^{\alpha,p} > 0$ \ such that for all
 \ $n \in \NN$,
 \[
   | R_{n,1,k}^{\alpha,p}(t) |
   \leq C_{1,k}^{\alpha,p} \frac{|t|^k}{n^{ (k - \alpha) / \alpha }}, \qquad
   t \in \RR .
 \]
Further, for arbitrary \ $\alpha \in (0,1) \cup (1,2)$ \ and \ $p \in (0,1)$,
 \ there exist \ $C_2^{\alpha,p} > 0$ \ and \ $C_3^{\alpha,p} > 0$ \ such that for
 all \ $n \in \NN$ \ and \ $|t| \leq C_2^{\alpha,p} n^{1/\alpha}$,
 \[
  \frac{| x_n^{\alpha,p}(t) |}{n} \leq \frac{1}{2} , \qquad
  | x_n^{\alpha,p}(t) | 
  \leq \begin{cases}
        C_3^{\alpha,p} \, |t|^\alpha , & \text{if \ $0 < \alpha < 1$,} \\[2mm]
        C_3^{\alpha,p} \, |t| \, n^{(\alpha - 1) / \alpha} ,
         & \text{if \ $1 < \alpha < 2$.}
       \end{cases}
 \]
\end{lem}

\begin{proof}
By the elementary inequality
 \begin{equation} \label{e^iu}
  \left| \ee^{\ii u} - \sum_{j=0}^{k-1} \frac{(\ii u)^j}{j!} \right|
  \leq \frac{|u|^k}{k!}, \qquad u\in\RR, 
 \end{equation}
 we obtain 
 \begin{align*}  
  | R_{n,1,k}^{\alpha,p}(t) |
  & \leq \frac{p \gamma_n}{q}
         \sum_{m = \lceil \log_r n \rceil}^\infty
          \frac{|t|^k}{k! \, r^{km/\alpha} \, \gamma_n^{k/\alpha}} r^m
    = \frac{p \, |t|^k}{k! \, q \, \gamma_n^{(k-\alpha)/\alpha}}
      \frac{q^{\lceil \log_r n \rceil (k-\alpha)/\alpha}}{1 - q^{(k-\alpha)/\alpha}} \\[2mm]
  & = \frac{p}{k! \, (q - q^{k / \alpha})}
      \frac{|t|^k}{(\gamma_n r^{\lceil \log_r n \rceil})^{(k-\alpha)/\alpha}}
    = - \frac{\mu_k^{\alpha,p}}{k!} \frac{|t|^k}{n^{(k-\alpha)/\alpha}},
 \end{align*}
 for all \ $t \in \RR$, \ hence the statement is satisfied with
 \ $C_{1,k}^{\alpha,p} = - \mu_k^{\alpha,p} / k! > 0$.

By \eqref{x}, Lemma \ref{y} and the estimate for \ $R_{n,1,2}^{\alpha,p}$,
 \begin{equation} \label{abs_x}
  | x_n^{\alpha,p}(t) |
  \leq C_1^{\alpha,p} |t|^\alpha
       + |\mu_1^{\alpha, p} | | t | n^{ ( \alpha - 1 ) / \alpha  }
       + C_{1,2}^{\alpha,p} \frac{|t|^2}{n^{ (2 - \alpha) / \alpha }} , \qquad
  t \in \RR ,
 \end{equation}
 hence for \ $|t| \leq C_2^{\alpha,p} n^{1/\alpha}$,
 \[
   \frac{| x_n^{\alpha,p}(t) |}{n}
   \leq C_1^{\alpha,p} ( C_2^{\alpha,p} )^\alpha
         + |\mu_1^{\alpha, p} | C_2^{\alpha,p}
         + C_{1,2}^{\alpha,p} ( C_2^{\alpha,p} )^2
   \leq \frac{1}{2}
 \]
 for all sufficiently small \ $C_2^{\alpha,p} > 0$.
\ If \ $0 < \alpha < 1$ \ then, by \eqref{abs_x}, for
 \ $|t| \leq C_2^{\alpha,p} n^{1/\alpha}$,
 \begin{align*}
  | x_n^{\alpha,p}(t) |
  & \leq |t|^\alpha
         \left\{ C_1^{\alpha,p}
                 + |\mu_1^{\alpha, p} | | t |^{1 - \alpha} n^{ ( \alpha - 1 ) / \alpha  }
                 + C_{1,2}^{\alpha,p} |t|^{2 - \alpha}
                   n^{ (\alpha - 2) / \alpha } \right\} \\[2mm]
  & \leq |t|^\alpha
         \left\{ C_1^{\alpha,p}
                 + |\mu_1^{\alpha, p} | ( C_2^{\alpha,p} )^{1 - \alpha}
                 + C_{1,2}^{\alpha,p} ( C_2^{\alpha,p} )^{2 - \alpha} \right\} .
 \end{align*}
If \ $1 < \alpha < 2$ \ then again by \eqref{abs_x}, for
 \ $|t| \leq C_2^{\alpha,p} n^{1/\alpha}$,
 \begin{align*}
  | x_n^{\alpha,p}(t) |
  & \leq |t| n^{ ( \alpha - 1 ) / \alpha  }
         \left\{ \frac{ C_1^{\alpha,p} |t|^{\alpha - 1} }{n^{ (\alpha - 1) / \alpha}}
                 + |\mu_1^{\alpha, p} |
                 + \frac{ C_{1,2}^{\alpha,p} |t| }{n^{ 1 / \alpha}} \right\}
    \\[2mm]
  & \leq |t| n^{ ( \alpha - 1 ) / \alpha  }
         \left\{ C_1^{\alpha,p} ( C_2^{\alpha,p} )^{\alpha - 1}
                 + |\mu_1^{\alpha, p} |
                 + C_{1,2}^{\alpha,p} C_2^{\alpha,p} \right\} ,
 \end{align*}
 establishing the case \ $1 < \alpha < 2$.
\end{proof}

Remark that Cs\"org\H{o} \cite{Cso_02} also derived these estimates for
 \ $|x_n^{\alpha,p}(t)|$.

For \ $k \in \{ 0,1,2,\dots \}$ \ and \ $|t| \leq C_2^{\alpha,p} n^{1/\alpha}$,
 \ introduce
 \begin{equation} \label{R2}
   R_{n,2,k}^{\alpha,p}(t)
   := \sum_{j=k}^\infty
       \frac{1}{j!}
       \big[ R_{n,1,2}^{\alpha,p}(t) + R_{n,3,2}^{\alpha,p}(t) \big]^j ,
 \end{equation}
 where
 \begin{equation} \label{R3}
   R_{n,3,k}^{\alpha,p}(t)
   := \sum_{j=k}^\infty
       \frac{(-1)^{j+1}}{j} \frac{[x_n^{\alpha,p}(t)]^j}{n^{j-1}} .
 \end{equation}
Then, by Lemma \ref{R1_x}, for \ $|t| \leq C_2^{\alpha,p} n^{1/\alpha}$, \ we may
 write
 \[
   n \log \left( 1 + \frac{x_n^{\alpha, p}(t)}{n} \right)
   = n \sum_{j=1}^\infty
        \frac{(-1)^{j+1}}{j} \left( \frac{x_n^{\alpha,p}(t)}{n} \right)^j
   = x_n^{\alpha,p}(t) + R_{n,3,2}^{\alpha,p}(t) .
 \]
Next we separate all the terms of order \ $n^{-\frac{k}{\alpha} - j}$ \ with
 \ $k + j \leq \ell - 1$ \ of the formal expansion of \ $\bff_n^{\alpha,p}$.
\ By \eqref{f_n^ap}, \eqref{x} and \eqref{R1}, for
 \ $|t| \leq C_2^{\alpha,p} n^{1/\alpha}$,
 \begin{align*}
  \bff_n^{\alpha,p}(t)
  & = \exp\left\{ - \ii t \mu_1^{\alpha, p} n^{ ( \alpha - 1 ) / \alpha  }
                  + x_n^{\alpha,p}(t)
                  + R_{n,3,2}^{\alpha,p}(t) \right\} \\[2mm]
  & = \exp\left\{ y_{\gamma_n}^{\alpha,p}(t) + R_{n,1,2}^{\alpha,p}(t)
                  + R_{n,3,2}^{\alpha,p}(t) \right\} \\[2mm]
  & = \ee^{y_{\gamma_n}^{\alpha,p}(t)}
      \sum_{j=0}^\infty
       \frac{1}{j!}
       \big[ R_{n,1,2}^{\alpha,p}(t) + R_{n,3,2}^{\alpha,p}(t) \big]^j \\[2mm]
  & = \ee^{y_{\gamma_n}^{\alpha,p}(t)}
      \Bigg\{ 1 + \sum_{j=1}^{\ell - 1}
                   \frac{1}{j!}
                   \big[ R_{n,1,2}^{\alpha,p}(t) + R_{n,3,2}^{\alpha,p}(t) \big]^j
              + R_{n,2,\ell}^{\alpha,p}(t) \Bigg\} \\[2mm]
  & = \ee^{y_{\gamma_n}^{\alpha,p}(t)}
      \Bigg\{ 1 + \sum_{j=1}^{\ell - 1}
                   \frac{1}{j!}
                   \Bigg[ \sum_{k=2}^{\ell-j+1}
                           \frac{\mu_k^{\alpha,p} \, (\ii t)^k }
                                {k! \, n^{ \frac{k}{\alpha} - 1 } } 
                          + \sum_{k=2}^{\ell-j+1}
                             \frac{(-1)^{k+1}}{k}
                             \frac{[x_n^{\alpha,p}(t)]^k}{n^{k-1}} + \\[2mm]
  & \phantom{= \ee^{y_{\gamma_n}^{\alpha,p}(t)}
               \Bigg\{ 1 + \Bigg[ \sum_{j=1}^{\ell - 1} \frac{1}{j!}}
                          + R_{n,1,\ell-j+2}^{\alpha,p}(t)
                          + R_{n,3,\ell-j+2}^{\alpha,p}(t) \Bigg]^j
                + R_{n,2,\ell}^{\alpha,p}(t) \Bigg\} .
 \end{align*}
Using notation
 \begin{equation} \label{tx}
   \tx_{n,m}^{\alpha,p}(t)
   := y_{\gamma_n}^{\alpha,p}(t)
      + \sum_{j=1}^m
         \frac{\mu_j^{\alpha,p}\,(\ii t)^j}
              {j!\,n^{ \frac{j}{\alpha} - 1 }} ,
   \qquad t \in \RR , \qquad m \in \NN,
 \end{equation}
 by \eqref{x}, we obtain
 \[
   x_n^{\alpha,p}(t)
   =  \tx_{n,m}^{\alpha,p}(t) + R_{n,1,m+1}^{\alpha,p}(t) ,
   \qquad t \in \RR , \qquad m \in \NN,
 \]
 and then
 \[
   \bff_n^{\alpha,p}(t)
   = \bg_{\alpha,p,\gamma_n}(t)
     \Bigg\{ 1 + \sum_{j=1}^{\ell - 1}
                  \frac{1}{j!}
                  \big[ \tR_{n,4,\ell-j+1}^{\alpha,p}(t) \big]^j
             + R_{n,5,\ell}^{\alpha,p}(t)
             + R_{n,2,\ell}^{\alpha,p}(t) \Bigg\} ,
 \]
 where, for \ $m \in\NN$ \ with \ $m \geq 2$,
 \begin{align}
  \tR_{n,4,m}^{\alpha,p}(t)
  & := \sum_{k=2}^m
        \frac{\mu_k^{\alpha,p} \, (\ii t)^k }
             {k! \, n^{ \frac{k}{\alpha} - 1 } }
       + \sum_{k=2}^m
          \frac{(-1)^{k+1}}{k}
          \frac{[\tx_{n,m-k+1}^{\alpha,p}(t)]^k}{n^{k-1}} , \label{tR4} \\[2mm]
  R_{n,5,\ell}^{\alpha,p}(t)
  & := \sum_{j=1}^{\ell - 1}
       \frac{1}{j!}
        \sum_{s=1}^j
         \binom{j}{s}
         \big[ R_{n,6,\ell-j+1}^{\alpha,p}(t) \big]^s \,
         \big[ \tR_{n,4,\ell-j+1}^{\alpha,p}(t) \big]^{j - s} , \label{R5}
 \end{align}
 with
 \begin{align}
  R_{n,6,m}^{\alpha,p}(t)
  & := R_{n,1,m+1}^{\alpha,p}(t) + \tR_{n,3,m}^{\alpha,p}(t)
       + R_{n,3,m+1}^{\alpha,p}(t) , \label{R6} \\[2mm]
  \tR_{n,3,m}^{\alpha,p}(t)
  & := \sum_{k=2}^m
        \frac{(-1)^{k+1}}{k \, n^{k-1}}
        \sum_{u=1}^k
         \binom{k}{u}
         \big[ R_{n,1,m-k+2}^{\alpha,p}(t) \big]^u \,
         \big[ \tx_{n,m-k+1}^{\alpha,p}(t) \big]^{k - u} . \label{tR3}
 \end{align}
Clearly \ $\sum_{j=1}^{\ell - 1} \big[ \tR_{n,4,\ell-j+1}^{\alpha,p}(t) \big]^j / j!$
 \ consists of certain terms of the formal infinite expansion
 \eqref{infinite_expantion} of \ $\bff_n^{\alpha,p}(t)$, \ and all the terms of
 order \ $n^{-\frac{k}{\alpha} - j}$ \ with \ $k + j \leq \ell - 1$ \ of
 \eqref{infinite_expantion} are contained in
 \ $\sum_{j=1}^{\ell - 1} \big[ \tR_{n,4,\ell-j+1}^{\alpha,p}(t) \big]^j / j!$. 
\ Consequently, for all \ $|t| \leq C_2^{\alpha,p} n^{1/\alpha}$,
 \begin{equation} \label{R7}
   \sum_{j=1}^{\ell - 1}
     \frac{1}{j!}
     \big[ \tR_{n,4,\ell-j+1}^{\alpha,p}(t) \big]^j 
    = \sum_{k=0}^{2 \ell-2} \,
        \sum_{j = - \lfloor k/2 \rfloor}^{\ell-k-1} \,
         \sum_{m = \max\{ 1, - j \}}^{k+j} \hspace*{-2mm}
          \frac{w_{m,k,j}^{\alpha, p} \, (- \ii t)^k \, [y_{\gamma_n}^{\alpha,p}(t)]^{j+m}}
               {m!  \, n^{\frac{k}{\alpha} + j}}
       + R_{n,7,\ell}^{\alpha,p}(t) , 
 \end{equation}
 where \ $R_{n,7,\ell}^{\alpha,p}(t)$ \ contains only terms of order
 \ $n^{-\frac{k}{\alpha} - j}$ \ with \ $k + j \geq \ell$.
\ Finally, we recognize that
 \begin{equation} \label{fg}
  \bff_n^{\alpha,p}(t)
  = \bg_{n,\ell-1}^{\alpha,p}(t)
    + \ee^{y_{\gamma_n}^{\alpha,p}(t)} R_{n,\ell}^{\alpha,p}(t) \qquad
    \text{for all \ $|t| \leq C_2^{\alpha,p} n^{1/\alpha}$,}
 \end{equation}
 where
 \[
   R_{n,\ell}^{\alpha,p}(t)
   := R_{n,7,\ell}^{\alpha,p}(t) + R_{n,5,\ell}^{\alpha,p}(t) + R_{n,2,\ell}^{\alpha,p}(t) .
 \]
In order to estimate the remainder term \ $R_{n,\ell}^{\alpha,p}$, \ we need the
 following lemmas.
Recall the definition \eqref{R3} of \ $R_{n,3,k}^{\alpha,p}$. 

\begin{lem} \label{R3_est}
For all \ $k \in \NN$, \ $\alpha \in (0,1) \cup (1,2)$, \ $p \in (0,1)$ \ and
 \ $n \in \NN$,
 \[
   | R_{n,3,k}^{\alpha,p}(t) |
   \leq \frac{2 |x_n^{\alpha,p}(t)|^k}{k \, n^{k - 1}}, \qquad
   \text{for all \ $|t| \leq C_2^{\alpha,p} n^{1/\alpha}$.}
 \]
\end{lem}

\begin{proof}
By Lemma \ref{R1_x},
 \[  
   | R_{n,3,k}^{\alpha,p}(t) |
   \leq \frac{|x_n^{\alpha,p}(t)|^k}{k \, n^{k-1}}
        \sum_{j=k}^\infty
         \left| \frac{x_n^{\alpha,p}(t)}{n} \right|^{j-k}
   \leq \frac{|x_n^{\alpha,p}(t)|^k}{k \, n^{k-1}}
        \sum_{m=0}^\infty
         \frac{1}{2^m} ,
 \]
 for all \ $|t| \leq C_2^{\alpha,p} n^{1/\alpha}$.
\end{proof}

\begin{lem} \label{R5_7}
For arbitrary \ $\ell \in \NN$, \ $\alpha \in (0,1) \cup (1,2)$ \ and
 \ $p \in (0,1)$, \ there exists \ $C_{2,\ell}^{\alpha,p} > 0$ \ such that for all
 \ $n \in \NN$ \ and \ $|t| \leq C_2^{\alpha,p} n^{1/\alpha}$,
 \[
   | R_{n,5,\ell}^{\alpha,p}(t) | + | R_{n,7,\ell}^{\alpha,p}(t) |
   \leq \begin{cases}
         \DS C_{2,\ell}^{\alpha,p}
             \frac{|t|^{(\ell+1)\alpha} + |t|^{2\ell\alpha}}{n^\ell} ,
          & \text{if \ $0 < \alpha < 1$,} \\[4mm]
         \DS C_{2,\ell}^{\alpha,p}
             \frac{|t|^{2 + (\ell-1)(2-\alpha)} + |t|^{2\ell}}{n^{\ell(2-\alpha)/\alpha}} ,
          & \text{if \ $1 < \alpha < 2$.}
        \end{cases}
 \]
\end{lem}

\begin{proof}
A term of order \ $n^{-\frac{k}{\alpha} - j}$ \ of \ $R_{n,7,\ell}^{\alpha,p}(t)$ \ is
 contained in the formal infinite expansion \eqref{infinite_expantion} of
 \ $\bff_n^{\alpha,p}(t)$, \ hence it has the form
 \[
   \frac{w_{m,k,j}^{\alpha, p} \, (- \ii t)^k \, [y_{\gamma_n}^{\alpha,p}(t)]^{j+m}}
        {m!  \, n^{\frac{k}{\alpha} + j}} ,
 \]
 where \ $k \geq 0$, \ $j \geq - \lfloor k/2 \rfloor$ \ and
 \ $1 \leq m \leq k+j$.
\ By Lemma \ref{y}, we have
 \[
   \left| \frac{(- \ii t)^k \, [y_{\gamma_n}^{\alpha,p}(t)]^{j+m}}
               {n^{\frac{k}{\alpha} + j}} \right|
   \leq (C_1^{\alpha,p})^{j+m}
        \frac{|t|^{k + (j+m)\alpha}}
             {n^{\frac{k}{\alpha} + j}} , \qquad n \in \NN, \quad t \in \RR .
 \]
For all \ $n \in \NN$ \ and \ $|t| \leq C_2^{\alpha,p} n^{1/\alpha}$,
 \[
   \frac{|t|^{k + (j+m)\alpha}}{n^{\frac{k}{\alpha} + j}}
   \leq \begin{cases}
         \DS (C_2^{\alpha,p})^{k(1-\alpha)}
             \frac{|t|^{(k+j+m)\alpha}}{n^{k+j}} ,
          & \text{if \ $0 < \alpha < 1$,} \\[4mm]
         \DS (C_2^{\alpha,p})^{(k+2j)(\alpha-1)}
             \frac{|t|^{(k+j)(2-\alpha)+m\alpha}}{n^{(k+j)(2-\alpha)/\alpha}} ,
          & \text{if \ $1 < \alpha < 2$.}
        \end{cases}
 \]
Since \ $R_{n,7,\ell}^{\alpha,p}(t)$ \ contains only terms of order
 \ $n^{-\frac{k}{\alpha} - j}$ \ with \ $k + j \geq \ell$ \ of the formal infinite
 expansion \eqref{infinite_expantion} of \ $\bff_n^{\alpha,p}(t)$, \ and
 \ $1 \leq m \leq k+j$, \ we obtain the estimate for
 \ $| R_{n,7,\ell}^{\alpha,p}(t) |$.
\ For the estimate of \ $| R_{n,5,\ell}^{\alpha,p}(t) |$ \ we derive from Lemmas
 \ref{y}, \ref{R1_x} and \ref{R3_est} the estimates
 \begin{align*}
  | \tx_{n,m}^{\alpha,p}(t) |
  & \leq \begin{cases}
          C_{3,m} |t|^\alpha, & \text{if  \ $0 < \alpha < 1$,} \\[2mm] 
          C_{3,m} |t| n^{(\alpha - 1) / \alpha},
           & \text{if  \ $1 < \alpha < 2$,}
         \end{cases} \\[2mm]
  | \tR_{n,3,m}^{\alpha,p}(t) |
  & \leq \begin{cases}
          \DS C_{3,m} \frac{|t|^{2 + (m - 1) \alpha}}{n^{\frac{2}{\alpha} + m - 2}},
           & \text{if  \ $0 < \alpha < 1$,} \\[4mm] 
          \DS C_{3,m} \frac{|t|^{m + 1}} {n^{\frac{m + 1}{\alpha} - 1}},
           & \text{if  \ $1 < \alpha < 2$,}
         \end{cases} \\[2mm]
  | \tR_{n,4,m}^{\alpha,p}(t) |
  & \leq \begin{cases}
          \DS C_{3,m} \frac{|t|^{2 \alpha}}{n},
           & \text{if  \ $0 < \alpha < 1$,} \\[4mm] 
          \DS C_{3,m} \frac{|t|^2} {n^{(2 - \alpha) / \alpha}},
           & \text{if  \ $1 < \alpha < 2$,}
         \end{cases} \\[2mm]
  | R_{n,6,m}^{\alpha,p}(t) |
  & \leq \begin{cases}
          \DS C_{3,m} \frac{|t|^{(m + 1) \alpha}}{n^m},
           & \text{if  \ $0 < \alpha < 1$,} \\[4mm] 
          \DS C_{3,m} \frac{|t|^{m + 1}} {n^{\frac{m + 1}{\alpha} - 1}},
           & \text{if  \ $1 < \alpha < 2$,}
         \end{cases}
 \end{align*}
 for all \ $n \in \NN$ \ and \ $|t| \leq C_2^{\alpha,p} n^{1/\alpha}$ \ and with
 sufficiently large constants \ $C_{3,m} > 0$.
\ The last two inequlaities imply the estimate for
 \ $| R_{n,5,\ell}^{\alpha,p}(t) |$.
\end{proof}

\begin{lem} \label{R2_est}
For arbitrary \ $\ell \in \NN$, \ $\alpha \in (0,1) \cup (1,2)$ \ and
 \ $p \in (0,1)$, \ there exist \ $C_{4,\ell}^{\alpha,p} > 0$ \ and
 \ $\vare_{\alpha,p} \in \big(0, C_2^{\alpha,p}\big]$ \ such that for all
 \ $n \in \NN$ \ and \ $|t| \leq \vare_{\alpha,p} \, n^{1/\alpha}$,
 \[
   |R_{n,2,\ell}^{\alpha,p}(t)|
   \leq \begin{cases}
         \DS C_{4,\ell}^{\alpha,p}
             \frac{|t|^{2\ell\alpha} \, \ee^{ C_1^{\alpha,p} |t|^\alpha / 2}}{n^\ell} ,
          & \text{if \ $0 < \alpha < 1$,} \\[4mm]
         \DS C_{4,\ell}^{\alpha,p}
             \frac{|t|^{2\ell} \, \ee^{ C_1^{\alpha,p} |t|^\alpha / 2}}
                  {n^{\ell(2-\alpha)/\alpha}} ,
         & \text{if \ $1 < \alpha < 2$.}
        \end{cases}
 \]
\end{lem}

\begin{proof}
Recalling the definition \eqref{R2} of \ $R_{n,2,\ell}^{\alpha,p}$, \ for all
 \ $n \in \NN$ \ and \ $|t| \leq C_2^{\alpha,p} n^{1/\alpha}$,
 \begin{align}
  | R_{n,2,\ell}^{\alpha,p}(t) |
  & \leq \sum_{j=\ell}^\infty
          \frac{1}{j!}
          \big( | R_{n,1,2}^{\alpha,p}(t) | + | R_{n,3,2}^{\alpha,p}(t) |\big)^j
    \nonumber \\[2mm]
  & \leq \frac{\big( | R_{n,1,2}^{\alpha,p}(t) | + | R_{n,3,2}^{\alpha,p}(t) |\big)^\ell}
              {\ell!}
         \sum_{j=\ell}^\infty
          \frac{\big( | R_{n,1,2}^{\alpha,p}(t) |
                      + | R_{n,3,2}^{\alpha,p}(t) |\big)^{j - \ell }}{(j - \ell)!} 
    \label{R2exp} \\[2mm]
  & = \frac{\big( | R_{n,1,2}^{\alpha,p}(t) | + | R_{n,3,2}^{\alpha,p}(t) |\big)^\ell}
              {\ell!}
         \ee^{ | R_{n,1,2}^{\alpha,p}(t) | + | R_{n,3,2}^{\alpha,p}(t) | }. \nonumber
 \end{align}
By Lemmas \ref{R1_x} and \ref{R3_est}, for all \ $n \in \NN$ \ and
 \ $|t| \leq C_2^{\alpha,p} n^{1/\alpha}$,
 \[
   | R_{n,1,2}^{\alpha,p}(t) | + | R_{n,3,2}^{\alpha,p}(t) |
   \leq C_{1,2}^{\alpha,p} \frac{t^2}{n^{(2-\alpha)/\alpha}}
        + \frac{|x_n^{\alpha,p}(t)|^2}{n} .
 \]
If \ $\alpha \in (0,1)$ \ then by Lemma \ref{R1_x}, for all \ $n \in \NN$ \ and
 \ $|t| \leq \vare_{\alpha,p} \, n^{1/\alpha}$ \ with
 \ $\vare_{\alpha,p} \in \big(0, C_2^{\alpha,p}\big]$,
 \begin{align}
  & | R_{n,1,2}^{\alpha,p}(t) | + | R_{n,3,2}^{\alpha,p}(t) |
    \leq C_{1,2}^{\alpha,p} \frac{t^2}{n^{(2-\alpha)/\alpha}}
         + \frac{(C_3^{\alpha,p})^2 \, |t|^{2 \alpha}}{n} \nonumber \\[2mm]
  & \leq \left[ C_{1,2}^{\alpha,p} (\vare_{\alpha,p})^{2-2\alpha}
                + (C_3^{\alpha,p})^2 \right]
         \frac{|t|^{2 \alpha}}{n} \label{R23_01_a} \\[2mm]
  & \leq \left[ C_{1,2}^{\alpha,p} (\vare_{\alpha,p})^{2-2\alpha}
                + (C_3^{\alpha,p})^2 \right]
         (\vare_{\alpha,p})^\alpha \, |t|^\alpha
    \leq \frac{1}{2} C_1^{\alpha,p} \, |t|^\alpha \label{R23_01_b}
 \end{align}
 for sufficiently small \ $\vare_{\alpha,p} > 0$.
\ Applying in the inequality \eqref{R2exp} the estimates \eqref{R23_01_a} and
 \eqref{R23_01_b} for
 \ $\big( | R_{n,1,2}^{\alpha,p}(t) | + | R_{n,3,2}^{\alpha,p}(t) |\big)^\ell$ \ and
 \ $\ee^{ | R_{n,1,2}^{\alpha,p}(t) | + | R_{n,3,2}^{\alpha,p}(t) | }$, \ respectively, we
 obtain the statement for \ $\alpha \in (0,1)$.

If \ $\alpha \in (1,2)$ \ then for all \ $n \in \NN$ \ and
 \ $|t| \leq \vare_{\alpha,p} \, n^{1/\alpha}$ \ with
 \ $\vare_{\alpha,p} \in \big(0, C_2^{\alpha,p}\big]$,
 \begin{align}
  | R_{n,1,2}^{\alpha,p}(t) | + | R_{n,3,2}^{\alpha,p}(t) |
  & \leq \left[ C_{1,2}^{\alpha,p} + (C_3^{\alpha,p})^2 \right]
         \frac{t^2}{n^{(2-\alpha)/\alpha}} \label{R23_12_a} \\[2mm]
  & \leq \left[ C_{1,2}^{\alpha,p} + (C_3^{\alpha,p})^2 \right]
         (\vare_{\alpha,p})^{2-\alpha} \, |t|^\alpha
    \leq \frac{1}{2} C_1^{\alpha,p} \, |t|^\alpha\label{R23_12_b}
 \end{align}
 for sufficiently small \ $\vare_{\alpha,p} > 0$. 
\ Applying in the inequality \eqref{R2exp} the estimates \eqref{R23_12_a} and
 \eqref{R23_12_b} for
 \ $\big( | R_{n,1,2}^{\alpha,p}(t) | + | R_{n,3,2}^{\alpha,p}(t) |\big)^\ell$ \ and
 \ $\ee^{ | R_{n,1,2}^{\alpha,p}(t) | + | R_{n,3,2}^{\alpha,p}(t) | }$, \ respectively, we
 obtain the statement for \ $\alpha \in (1,2)$.
\end{proof}

By \eqref{fg}, the first inequality of Lemma \ref{y}, and Lemmas \ref{R5_7} and
 \ref{R2_est}, we obtain

\begin{lem} \label{R}
For arbitrary \ $\ell \in \NN$, \ $\alpha \in (0,1) \cup (1,2)$ \ and
 \ $p \in (0,1)$, \ there exists \ $C_{5,\ell}^{\alpha,p} > 0$ \ such that for all
 \ $n \in \NN$ \ and \ $|t| \leq \vare_{\alpha,p} \, n^{1/\alpha}$,
 \[
   | \bff_n^{\alpha,p}(t) - \bg_{n,\ell - 1}^{\alpha,p}(t) |
   \leq \begin{cases}
         \DS C_{5,\ell}^{\alpha,p}
             \frac{|t|^{(\ell+1)\alpha} + |t|^{2\ell\alpha}}{n^\ell}
             \ee^{ - C_1^{\alpha,p} |t|^\alpha / 2} ,
          & \text{if \ $\alpha \in (0,1)$,} \\[4mm]
         \DS C_{5,\ell}^{\alpha,p}
             \frac{|t|^{2 + (\ell-1)(2-\alpha)} + |t|^{2\ell}}{n^{\ell(2-\alpha)/\alpha}}
             \ee^{ - C_1^{\alpha,p} |t|^\alpha / 2} ,
         & \text{if \ $\alpha \in (1,2)$.}
        \end{cases}
 \]
\end{lem}

Now by Lemma \ref{Esseen1},
 \[
   \Delta_{n,\ell}^{\alpha,p}
   := \sup_{ x \in \RR }
       \Bigg| \PP \bigg\{ \frac{S_n - \mu_1^{ \alpha, p } n }
                               { n^{ 1 / \alpha } } \leq x \bigg\}
              - G_{n,\ell - 1}^{\alpha, p}(x) \Bigg|
   \leq \frac{b}{\pi} \Delta_{n,\ell,1}^{\alpha,p} + c_b \Delta_{n,\ell,2}^{\alpha,p}
 \]
 for every \ $b > 1$, \ where
 \begin{align*}
  \Delta_{n,\ell,1}^{\alpha,p}
  & := \int_0^{\vare_{\alpha,p} \, n^{1 / \alpha}}
        \frac{| \bff_n^{\alpha,p}(t) - \bg_{n,\ell - 1}^{\alpha,p}(t) |}{|t|} \, \dd t,
    \\[2mm]
  \Delta_{n,\ell,2}^{\alpha,p}
  & := \frac{M_{n,\ell}^{\alpha,p}}{\vare_{\alpha,p} \, n^{1 / \alpha}}
  \qquad \text{with} \quad
  M_{n,\ell}^{\alpha,p}
  := \sup_{x \in \RR} \left| \frac{\dd G_{n,\ell - 1}^{\alpha,p} (x)}{\dd x} \right| .
 \end{align*}
By Lemmas 4 and 6 of Cs\"org\H{o} \cite{Cso_07b},
 \ $\sup_{n \in \NN} M_{n,\ell}^{\alpha,p} < \infty$, \ and hence
 \ $\Delta_{n,\ell,2}^{\alpha,p} = O\left( \frac{1}{n^{1 / \alpha}} \right)$ \ for all
 fixed \ $\ell \in \NN$.
\ Using the simple fact that
 \begin{equation} \label{int1}
  \int_0^\infty t^\beta \ee^{-C t^\alpha} \, \dd t < \infty , \qquad
  \beta > -1 , \qquad C > 0 ,
 \end{equation}
 by Lemma \ref{R} we obtain
 \begin{equation} \label{Delta1}
   \Delta_{n,\ell,1}^{\alpha,p}
   = O \left( \frac{1}{n^\ell} + \frac{1}{n^{\ell(2-\alpha)/\alpha}} \right) .
 \end{equation}
Consequently, we conclude
 \ $\Delta_{n,\ell}^{\alpha,p}
    = O \left( \frac{1}{n^\ell} + \frac{1}{n^{\ell(2-\alpha)/\alpha}}
               + \frac{1}{n^{1 / \alpha}} \right)$,
 \ which implies the first statement.

The reduction of the order \ $O\left(\frac{1}{n^{1/\alpha}}\right)$ \ to 
 \ $o\left(\frac{1}{n^{1/\alpha}}\right)$ \ when \ $\ell \in \{2,3,\dots\}$,
 \ $\alpha \in \left( \frac{1}{\ell} , 1 \right)
               \cup \left(1, 2 - \frac{1}{\ell} \right) $
 \ and \ $r^{1/\alpha} \notin \NN$ \ is based on the following classical result
 due to Esseen \cite{Esseen_45}.

\begin{lem}[Esseen] \label{Esseen2}
If \ $\bff$ \ is the characteristic function of a non-lattice distribution,
 then for every fixed \ $\vare>0$, \ there exists a sequence
 \ $\lambda_n\to\infty$ \ such that
 \[
   \int_\vare^{\lambda_n} \frac{|\bff(t)|^n}{t} \, \dd t
   = o \left( \ee^{-\sqrt{n}/2} \right)
   \qquad \text{as \ $n\to\infty$.}
 \]
\end{lem}

It is easy to check that the distribution of \ $X$ \ is non-lattice if and only
 if \ $r^{1/\alpha} \notin \NN$. 
\ Thus there exists a sequence \ $\lambda_n^{\alpha,p}\to\infty$ \ such that
 \begin{align*}
  I_n^{\alpha,p}
  &:= \int_{\vare_{\alpha,p} \, n^{1 / \alpha}}^{\lambda_n^{\alpha,p} \, n^{1/\alpha}}
      \frac{|\bff_n^{\alpha,p}(t)|}{t} \, \dd t
   = \int_{\vare_{\alpha,p} \, n^{1 / \alpha}}^{\lambda_n^{\alpha,p} \, n^{1/\alpha}}
      \frac{|\bff_{\alpha,p}( t/ n^{1/\alpha} )|^n}{t} \, \dd t \\[2mm]
  &= \int_{\vare_{\alpha,p}}^{\lambda_n^{\alpha,p}}
      \frac{|\bff_{\alpha,p}(s)|^n}{s} \, \dd s
   = o \left( \ee^{-\sqrt{n}/2} \right)
   = o \left( \frac{1}{n^{1/\alpha}} \right) \qquad
     \text{as \ $n\to\infty$.}
 \end{align*}
Clearly \eqref{Delta1} implies
 \ $\Delta_{n,\ell,1}^{\alpha,p} = o \left( \frac{1}{n^{1/\alpha}} \right)$ \ as
 \ $n\to\infty$ \ for
 \ $\alpha \in \left( \frac{1}{\ell} , 1 \right)
               \cup \left(1, 2 - \frac{1}{\ell} \right)$.
\ Hence by Lemma \ref{Esseen1} now with \ $T = \lambda_n^{\alpha,p} \, n^{1/\alpha}$
 \ we have
 \[
   \Delta_{n,\ell}^{\alpha,p}
   \leq \frac{b}{\pi} 
        \left( \Delta_{n,\ell,1}^{\alpha,p}
               + I_n^{\alpha,p}
               + I_{n,\ell}^{\alpha,p} \right)
        + c_b \frac{M_{n,\ell}^{\alpha,p}}{\lambda_n^{\alpha,p} n^{1/\alpha}}
   = \frac{b}{\pi} I_{n,\ell}^{\alpha,p} + o \left( \frac{1}{n^{1/\alpha}} \right)
 \]
 as \ $n\to\infty$, \ where
 \[
   I_{n,\ell}^{\alpha,p}
   := \int_{\vare_{\alpha,p} \, n^{1 / \alpha}}^\infty
       \frac{|\bg_{n,\ell-1}^{\alpha,p}(t)|}{t} \, \dd t
   = O \left( n^{\ell - 2} \ee^{ - C_1^{\alpha,p} (\vare_{\alpha,p})^\alpha \, n} \right)
   = o \left( \frac{1}{n^{1/\alpha}} \right)
 \]
 as \ $n\to\infty$, \ using Lemma \ref{y} and the simple fact that
 \begin{equation} \label{int2}
   \int_n^\infty u^\beta \ee^{- u} \, \dd u = O(n^\beta \ee^{- n}) , \qquad
   \beta \in \RR .
 \end{equation}
Thus \ $\Delta_{n,\ell}^{\alpha,p} = o \left( \frac{1}{n^{1/\alpha}} \right)$ \ as
 \ $n\to\infty$ \ for
 \ $\alpha \in \left( \frac{1}{\ell} , 1 \right)
               \cup \left(1, 2 - \frac{1}{\ell} \right)$.

\section{Nonuniform bounds in asymptotic expansions}
\label{nonunif_bounds}

The main results are contained in the following

\begin{pro} \label{main_nonuniform}
For \ $\ell \in \{2,3,\dots\}$ \ and
 \ $\alpha \in \left( 1 , 2 - \frac{1}{\ell} \right]$,
 \[
   \sup_{ x \in \RR } \,
    ( 1 + |x| ) 
   \Bigg| \PP \bigg\{ \frac{S_n - \mu_1^{ \alpha, p } n } { n^{ 1 / \alpha } }
                      \leq x \bigg\}
          - G_{n,\ell - 1}^{\alpha, p}(x) \Bigg|
   = O \left( \frac{1}{n^{1 / \alpha}} \right) .
 \] 
For \ $\ell \in \NN$ \ and
 \ $\alpha \in \left( 2 - \frac{1}{\ell} , \, 2 \right)$,
 \[
   \sup_{ x \in \RR } \,
    ( 1 + |x| ) 
   \Bigg| \PP \bigg\{ \frac{S_n - \mu_1^{ \alpha, p } n } { n^{ 1 / \alpha } }
                      \leq x \bigg\}
          - G_{n,\ell - 1}^{\alpha, p}(x) \Bigg|
   = O \left( \frac{1}{ n^{ \ell ( 2 - \alpha ) / \alpha } } \right) .
 \] 
For \ $\ell \in \{2,3,\dots\}$,
 \ $\alpha \in \left( 1, \, 2 - \frac{1}{\ell} \right)$ \ and
 \ $r^{1/\alpha} \notin \NN$,
 \[
   \sup_{ x \in \RR } \,
    ( 1 + |x| ) 
   \Bigg| \PP \bigg\{ \frac{S_n - \mu_1^{ \alpha, p } n } { n^{ 1 / \alpha } }
                      \leq x \bigg\}
          - G_{n,\ell - 1}^{\alpha, p}(x) \Bigg|
   = o \left( \frac{1}{ n^{ \ell ( 2 - \alpha ) / \alpha } } \right) .
 \] 
\end{pro}

For \ $\ell = 1$ \ and \ $\alpha \in (1,2)$, \ we have a nonuniform rate of
 merge
 \[
   \sup_{ x \in \RR } \,
    ( 1 + |x| ) 
   \Bigg| \PP \bigg\{ \frac{S_n - \mu_1^{ \alpha, p } n } { n^{ 1 / \alpha } }
                      \leq x \bigg\}
          - G_{\alpha, p,\gamma_n}(x) \Bigg|
   = O \left( \frac{1}{ n^{ ( 2 - \alpha ) / \alpha } } \right) .
 \]
Certain terms of \ $G_{n,\ell - 1}^{\alpha, p}$ \ are of the same or of a smaller
 order than the remainder terms \ $O \left( \frac{1}{n^{1/\alpha}} \right)$ \ or
 \ $o \left( \frac{1}{n^{1/\alpha}} \right)$.
\ Using the boundedness of the functions \ $G_{\alpha, p, \gamma_n}^{(k,j)}$, \ the
 expansions may be simplified as follows.
Recall the definition of the approximation functions \ $\tG_{n,\ell}^{\alpha, p}$
 \ and \ $\ttG_{n,\ell}^{\alpha, p}$ \ from Theorem \ref{MAIN_uniform}.
 
\begin{thm} \label{MAIN_nonuniform} 
For \ $\ell \in \{ 2, 3, \dots\}$ \ and
 \ $ \alpha \in \left( 2 - \frac{1}{\ell-1} , \, 2 - \frac{1}{\ell} \right] $,
 \[
   \sup_{ x \in \RR } \,
    ( 1 + |x| ) 
    \Bigg| \PP \bigg\{ \frac{S_n - c_n^{ \alpha, p } } { n^{ 1 / \alpha } }
                       \leq x \bigg\}
           - \tG_{n,\ell - 1}^{\alpha, p}(x) \Bigg|
   = O \left( \frac{1}{n^{1 / \alpha}} \right) .
 \]
For \ $\ell = 2$, \ $\alpha \in \left( 1, \frac{3}{2} \right)$, \ or for
 \ $\ell \in \{ 3, 4, \dots\}$,
 \ $ \alpha \in \left[ 2 - \frac{1}{\ell-1} , \, 2 - \frac{1}{\ell} \right)$,
 \ and for \ $ r^{ 1 / \alpha} \notin \NN $,
 \[
   \sup_{ x \in \RR } \,
    ( 1 + |x| ) 
    \Bigg| \PP \bigg\{ \frac{S_n \! - \mu_1^{ \alpha, p } n }
                            { n^{ 1 / \alpha } } \!
                       \leq x \bigg\}
           - \ttG_{n,\ell - 1}^{\alpha, p}(x) \Bigg|
   = o \left( \frac{1}{n^{1 / \alpha}} \right) ,
 \]
\end{thm}

\section{Proof of Proposition \ref{main_nonuniform}}

Fix \ $\ell \in \NN$ \ and \ $\alpha \in (1,2)$.
\ The proof is based on a result due to Osipov \cite{Osi_72} (see Petrov
 \cite[Lemmas 6.7, 6.8]{Pet_75} and Hall \cite[Lemma 1.1]{Hall_82}; note that
 the result has been stated only for \ $k\geq2$, \ but the proof works for
 \ $k=1$ \ as well).

\begin{lem}[Osipov] \label{Osipov}
Let \ $F$ \ be a distribution function and \ $G$ \ be a differentiable function
 of bounded variation on \ $\RR$ \ with Fourier--Stieltjes transforms
 \ $ \bff(t) = \int_{-\infty}^\infty \ee^{\ii t x} \, \dd F(x)$ \ and
 \ $ \bg(t) = \int_{-\infty}^\infty \ee^{\ii t x} \, \dd G(x)$, \ $t \in\RR$, \ such
 that \ $G(-\infty) = \lim_{x\to-\infty} G(x) = 0$,
 \ $G(+\infty) = \lim_{x\to+\infty} G(x) = 1$, \ and
 \[
   \int_{-\infty}^\infty |x|^k \, | \dd (F(x)-G(x)) | < \infty 
 \]
 with some \ $k \in \NN$.
\ Then the function \ $H_k(x) := x^k (F(x)-G(x))$, \ $x \in \RR$, \ is of
 bounded variation on \ $\RR$ \ with Fourier--Stieltjes transform
 \[
   \bh_k(t)
   = \int_{-\infty}^\infty \ee^{\ii t x} \, \dd H_k(x)
   = \frac{k!}{(-\ii t)^k}
     \sum_{j=0}^k
      \frac{(-t)^j}{j!} \frac{\dd^j(\bff-\bg)(t)}{\dd t^j} ,
 \]
 for \ $t \in \RR$ \ with \ $t \neq 0$, \ and there exists \ $c_k > 0$ \ such
 that
 \begin{align*}
  & \sup_{x \in \RR} \, ( 1 + |x|^k ) | F(x) - G(x) | \\[2mm]
  & \leq c_k \int_{-T}^T \left| \frac{\bff(t) - \bg(t)}{t} \right| \, \dd t
         + c_k \int_{-T}^T \left| \frac{\bh_k(t)}{t} \right| \, \dd t
         + c_k \frac{\sup_{x \in \RR} ( 1 + |x|^k ) | G'(x) |}{T}         
 \end{align*}
 for every choice of \ $T > 0$.
\end{lem}

We apply this lemma for
 \ $F(x) = F_n^{\alpha, p}(x)
         = \PP \left\{ \frac{S_n - \mu_1^{ \alpha, p } n } { n^{ 1 / \alpha } }
                       \leq x \right\}$
 \ and \ $G(x) = G_{n,\ell - 1}^{\alpha, p}(x)$, \ $x \in \RR$, \ introduced in
 Section \ref{unif_bounds}.
By Lemma 6 of Cs\"org\H{o} \cite{Cso_07b}, \ $G_{n,\ell - 1}^{\alpha, p}$ \ is a
 differentiable function of bounded variation on \ $\RR$ \ with
 \ $G_{n,\ell - 1}^{\alpha, p}(-\infty) = 0$ \ and
 \ $G_{n,\ell - 1}^{\alpha, p}(+\infty) = 1$.
\ Since \ $\alpha \in (1,2)$, \ the expectation \ $\EE(X)$ \ of the gain in one
 game is finite, which implies
 \ $\int_{-\infty}^\infty |x| \, \dd F_n^{\alpha, p}(x)
    = \EE \left|\frac{S_n - \mu_1^{ \alpha, p } n } { n^{ 1 / \alpha } }\right|
    < \infty $.
\ By Lemma 6 of Cs\"org\H{o} \cite{Cso_07b}, \ $\alpha \in (1,2)$ \ also
 implies
 \ $\int_{-\infty}^\infty |x| \, | \dd G_{n,\ell - 1}^{\alpha, p}(x) | < \infty$,
 \ since \ $G_{n,\ell - 1}^{\alpha, p}$ \ is a linear combination of the continuously
 differentiable functions
 \ $G_{\alpha, p, \gamma}^{(k,j)}$, \ $k,j\in\{0,1,2,\dots\}$, \ $\gamma \in (q,1]$,
 \ and
 \ $\int_{-\infty}^\infty |x| \, | \dd G_{\alpha, p, \gamma}^{(k,j)}(x) |
    = \int_{-\infty}^\infty |x| | G_{\alpha, p, \gamma}^{(k+1,j)}(x) | \, \dd x < \infty$.
\ Consequently,
 \[
  \int_{-\infty}^\infty |x| \, | \dd (F_n^{\alpha, p}(x)-G_{n,\ell - 1}^{\alpha, p}(x)) | \\
  \leq \int_{-\infty}^\infty |x| \, \dd F_n^{\alpha, p}(x)
       + \int_{-\infty}^\infty |x| \, | G_{n,\ell - 1}^{\alpha, p}(x)) |
  < \infty .
 \]
By Lemma \ref{Osipov}, the function
 \ $H_{n,\ell - 1}^{\alpha, p}(x) := x (F_n^{\alpha, p}(x)-G_{n,\ell - 1}^{\alpha, p}(x))$,
 \ $x \in \RR$, \ is of bounded variation on \ $\RR$ \ with Fourier--Stieltjes
 transform
 \[
   \bh_{n,\ell - 1}^{\alpha, p}(t)
   = \int_{-\infty}^\infty \ee^{\ii t x} \, \dd H_{n,\ell - 1}^{\alpha, p}(x)
   = \frac{\ii}{t}(\bff_n^{\alpha, p}(t)-\bg_{n,\ell - 1}^{\alpha, p}(t))
     - \ii \frac{\dd(\bff_n^{\alpha, p}-\bg_{n,\ell - 1}^{\alpha, p})(t)}{\dd t}
 \]
 for \ $t \in \RR$ \ with \ $t \neq 0$. 
\ Now we apply Lemma \ref{Osipov} for \ $T = \vare_{\alpha,p} n^{1/\alpha}$ \ with
 \ $\vare_{\alpha,p} > 0$ \ from Lemmas \ref{R2_est} and \ref{R}, and we obtain
 \begin{align*}
   \tDelta_{n,\ell}^{\alpha,p}
   & := \sup_{ x \in \RR }
         ( 1 + |x| )
         \Bigg| \PP \bigg\{ \frac{S_n - \mu_1^{ \alpha, p } n }
                                 { n^{ 1 / \alpha } } \leq x \bigg\}
                - G_{n,\ell - 1}^{\alpha, p}(x) \Bigg| \\[2mm]
   & \leq 2 c_1 \Big( \tDelta_{n,\ell,1}^{\alpha,p}
                      + \tDelta_{n,\ell,2}^{\alpha,p}
                      + \tDelta_{n,\ell,3}^{\alpha,p} \Big)
          + c_1 \tDelta_{n,\ell,4}^{\alpha,p} ,
 \end{align*}
 where
 \begin{align*}
  \tDelta_{n,\ell,1}^{\alpha,p}
  & := \int_0^{\vare_{\alpha,p} \, n^{1 / \alpha}}
        \frac{| \bff_n^{\alpha,p}(t) - \bg_{n,\ell - 1}^{\alpha,p}(t) |}{|t|} \, \dd t,
    \\[2mm]
  \tDelta_{n,\ell,2}^{\alpha,p}
  & := \int_0^{\vare_{\alpha,p} \, n^{1 / \alpha}}
        \frac{| \bff_n^{\alpha,p}(t) - \bg_{n,\ell - 1}^{\alpha,p}(t) |}{|t|^2} \, \dd t,
    \\[2mm]
  \tDelta_{n,\ell,3}^{\alpha,p}
  & := \int_0^{\vare_{\alpha,p} \, n^{1 / \alpha}}
        \frac{1}{|t|}
        \left| \frac{\dd(\bff_n^{\alpha, p}-\bg_{n,\ell - 1}^{\alpha, p})(t)}{\dd t}
        \right| \, \dd t,
    \\[2mm]
  \tDelta_{n,\ell,4}^{\alpha,p}
  & := \frac{\tM_{n,\ell}^{\alpha,p}}{\vare_{\alpha,p} \, n^{1 / \alpha}}
  \qquad \text{with} \quad
  \tM_{n,\ell}^{\alpha,p}
  := \sup_{x \in \RR} \,
      ( 1 + |x| ) \left| \frac{\dd G_{n,\ell - 1}^{\alpha,p} (x)}{\dd x} \right| .
 \end{align*}
By Lemmas 4 and 6 of Cs\"org\H{o} \cite{Cso_07b},
 \ $\sup_{n \in \NN} \tM_{n,\ell}^{\alpha,p} < \infty$, \ and hence
 \ $\tDelta_{n,\ell,4}^{\alpha,p} = O\left( \frac{1}{n^{1 / \alpha}} \right)$ \ for all
 fixed \ $\ell \in \NN$.
\ By \eqref{Delta1}, we have
 \ $\tDelta_{n,\ell,1}^{\alpha,p} = O \left( \frac{1}{n^{\ell(2-\alpha)/\alpha}} \right)$.
\ By \eqref{fg} and Lemmas \ref{y}, \ref{R5_7} and \ref{R2_est} we obtain
 \[
   \tDelta_{n,\ell,2}^{\alpha,p}
   \leq \int_0^{\vare_{\alpha,p} \, n^{1 / \alpha}}
         \frac{ | R_{n,\ell}^{\alpha,p}(t) | \, \ee^{- C_1^{\alpha,p} |t|^\alpha } }
              {|t|^2} \,
         \dd t
   = O \left( \frac{1}{n^{\ell(2-\alpha)/\alpha}} \right) . 
 \]
The aim of the following discussion is to find an appropriate estimate for
 \ $\left| \frac{\dd(\bff_n^{\alpha, p}-\bg_{n,\ell - 1}^{\alpha, p})(t)}{\dd t}
    \right|$
 \ if \ $|t| \leq \vare_{\alpha,p} \, n^{1 / \alpha}$.
\ Using formula \eqref{fg}, first we calculate the derivatives of the
 ingredients of \ $\bff_n^{\alpha, p}-\bg_{n,\ell - 1}^{\alpha, p}$. 
\ The function \ $y_\gamma^{\alpha,p}$ \ is differentiable and
 \begin{equation} \label{der_y}
  \frac{\dd y_\gamma^{\alpha,p}(t)}{\dd t}
  = \sum\limits_{k = - \infty}^\infty
     \Bigg( \exp \left\{ \frac{\ii t r^{k / \alpha}}{\gamma^{1 / \alpha}} \right\}
            -1 \Bigg)
     \frac{\ii p \gamma^{(\alpha-1)/\alpha}}{q r^{k(\alpha-1)/\alpha}} , \qquad
  t \in \RR,
 \end{equation}
 since \ $\alpha \in (1,2)$ \ implies absolute convergence of this series.
Differentiability of \ $y_\gamma^{\alpha,p}$ \ implies that for each \ $m \in \NN$,
 \ the function \ $\tx_{n,m}^{\alpha,p}$, \ introduced in \eqref{tx}, hence, for
 each \ $m \geq 2$, \ the function \ $\tR_{n,4,m}^{\alpha,p}$, \ given in
 \eqref{tR4}, and hence for each \ $k \in \NN$, \ the function
 \ $R_{n,7,k}^{\alpha,p}$, \ defined in \eqref{R7}, is differentiable on the whole
 \ $\RR$.
\ For each \ $k \geq 2$, \ the function \ $R_{n,1,k}^{\alpha,p}$, \ given in
 \eqref{R1k}, is differentiable and
 \begin{equation} \label{der_R1}
   \frac{\dd R_{n,1,k}^{\alpha,p}(t)}{\dd t}
   = - \sum_{m = \lceil \log_r n \rceil}^\infty 
         \Bigg( \exp \left\{ \frac{\ii t}{r^{m/\alpha} \gamma_n^{1/\alpha}} \right\}
                 - \sum_{j=0}^{k-2}
                    \frac{(\ii t)^j}
                         {j! \, r^{jm/\alpha} \, \gamma_n^{j/\alpha}} \Bigg)
          \frac{\ii p \gamma_n^{(\alpha-1)/\alpha}}{q} r^{m(\alpha-1)/\alpha} , 
 \end{equation}
 for \ $t \in \RR$, \ since \ $\alpha \in (1,2)$ \ implies absolute
 convergence of this series.
Differentiability of \ $\tx_{n,m}^{\alpha,p}$, \ $m \in \NN$, \ and
 \ $R_{n,1,k}^{\alpha,p}$, \ $k \geq 2$, \ imply that for each \ $m \geq 2$,
 \ the function \ $\tR_{n,3,m}^{\alpha,p}$, \ given in \eqref{tR3}, is
 differentiable on the whole \ $\RR$.
\ Recalling formula \eqref{x}, we obtain differentiability of the function
 \ $x_n^{\alpha, p}$, \ and
 \begin{equation} \label{der_x}
   \frac{\dd x_n^{\alpha, p}(t)}{\dd t}
   = \frac{\dd y_{\gamma_n}^{\alpha,p}(t)}{\dd t}
     + \ii \mu_1^{\alpha, p} n^{( \alpha - 1 ) / \alpha} 
     + \frac{\dd R_{n,1,2}^{\alpha,p}(t)}{\dd t} , \qquad t \in \RR .
 \end{equation}
For each \ $k \in \NN$, \ the function \ $R_{n,3,k}^{\alpha,p}$, \ defined in
 \eqref{R3}, is differentiable for \ $|t| < C_2^{\alpha,p} \, n^{1/\alpha}$ \ and
 \begin{equation} \label{der_R3}
   \frac{\dd R_{n,3,k}^{\alpha,p}(t)}{\dd t}
   = \frac{\dd x_n^{\alpha, p}(t)}{\dd t}
     \sum_{j=k-1}^\infty
       \frac{(-1)^j [x_n^{\alpha,p}(t)]^j}{n^j} , \qquad
   |t| < C_2^{\alpha,p} \, n^{1/\alpha},
 \end{equation}
 since by Lemma \ref{R1_x}, \ $\alpha \in (1,2)$ \ implies absolute convergence
 of this series for \ $|t| < C_2^{\alpha,p} \, n^{1/\alpha}$.
\ For each \ $k \in \NN$, \ the function \ $R_{n,2,k}^{\alpha,p}$, \ given in
 \eqref{R2}, is differentiable for \ $|t| < C_2^{\alpha,p} \, n^{1/\alpha}$ \ and
 \begin{equation} \label{derR2}
  \frac{\dd R_{n,2,k}^{\alpha,p}(t)}{\dd t}
  = \left( \frac{\dd R_{n,1,2}^{\alpha,p}(t)}{\dd t}
           + \frac{\dd R_{n,3,2}^{\alpha,p}(t)}{\dd t} \right)
  R_{n,2,k-1}^{\alpha,p}(t) ,
 \end{equation}
 for \ $|t| < C_2^{\alpha,p} \, n^{1/\alpha}$.
\ Differentiability of \ $R_{n,1,k}^{\alpha,p}$, \ $k \geq 2$,
 \ $\tR_{n,3,m}^{\alpha,p}$, \ $m \in \NN$, \ and
 \ $R_{n,3,k}^{\alpha,p}$, \ $k \in \NN$, \ imply that for each \ $m \geq 2$,
 \ the function \ $R_{n,6,m}^{\alpha,p}$, \ given in \eqref{R6}, and hence for each
 \ $k \in \NN$, \ the function \ $R_{n,5,k}^{\alpha,p}$, \ introduced in
 \eqref{R5}, is differentiable for \ $|t| < C_2^{\alpha,p} \, n^{1/\alpha}$.
\ Consequently by \eqref{fg}, we conclude
 \begin{align*} \label{der_f-g}
  \frac{\dd(\bff_n^{\alpha, p}-\bg_{n,\ell - 1}^{\alpha, p})(t)}{\dd t}
  & = \ee^{y_{\gamma_n}^{\alpha,p}(t)}
      \left( \frac{\dd R_{n,2,\ell}^{\alpha,p}(t)}{\dd t}
             + \frac{\dd R_{n,5,\ell}^{\alpha,p}(t)}{\dd t}
             + \frac{\dd R_{n,7,\ell}^{\alpha,p}(t)}{\dd t} \right) \\[2mm]
  & \quad + \frac{\dd y_{\gamma_n}^{\alpha,p}(t)}{\dd t}
            \big( \bff_n^{\alpha, p}(t)-\bg_{n,\ell - 1}^{\alpha, p}(t) \big) , \qquad
  |t| < C_2^{\alpha,p} \, n^{1/\alpha} .
 \end{align*}
In order to estimate
 \ $\left| \frac{\dd(\bff_n^{\alpha, p}-\bg_{n,\ell - 1}^{\alpha, p})(t)}{\dd t}
    \right|$,
 \ we need the following lemmas.

\begin{lem} \label{der_y_est}
For arbitrary \ $\alpha \in (1,2)$ \ and \ $p \in (0,1)$, \ there exists
 \ $C_4^{\alpha,p} > 0$ \ such that, uniformly in \ $\gamma \in (q,1]$,
 \[
   \Im \left( \frac{\dd y_\gamma^{\alpha,p}(t)}{\dd t} \right)
   \leq - C_4^{\alpha,p} \, |t|^{\alpha-1}, \qquad
   \left| \frac{\dd y_\gamma^{\alpha,p}(t)}{\dd t} \right|
   \leq C_4^{\alpha,p} \, |t|^{\alpha-1}, \qquad  t \in \RR .
 \]
\end{lem}

\begin{proof}
These estimates can be derived in the same way as the inequality (4) in the
 proof of Lemma 3 of Cs\"org\H{o} \cite{Cso_02} and the second statement of
 Lemma 3 in Cs\"org\H{o} \cite{Cso_07b}, using \eqref{der_y}.
\end{proof}

The following lemmas can be proved as Lemmas \ref{R1_x}, \ref{R3_est},
 \ref{R5_7}, \ref{R2_est} and \ref{R}, respectively.

\begin{lem} \label{der_R1_x_est}
For arbitrary \ $k \geq 2$, \ $\alpha \in (1,2)$ \ and \ $p \in (0,1)$, \ there
 exists \ $C_{5,k}^{\alpha,p} > 0$ \ such that for all \ $n \in \NN$,
 \[
   \left| \frac{\dd R_{n,1,k}^{\alpha,p}(t)}{\dd t} \right|
   \leq C_{5,k}^{\alpha,p} \frac{|t|^{k-1}}{n^{ (k - \alpha) / \alpha }}, \qquad
   t \in \RR .
 \]
Further, for arbitrary \ $\alpha \in (1,2)$ \ and \ $p \in (0,1)$, \ there
 exists \ $C_5^{\alpha,p} > 0$ \ such that for all \ $n \in \NN$ \ and
 \ $|t| \leq C_2^{\alpha,p} \, n^{1/\alpha}$,
 \[
   \left| \frac{\dd x_n^{\alpha,p}(t)}{\dd t} \right|
   \leq C_5^{\alpha,p} \, n^{(\alpha - 1) / \alpha} .
 \]
\end{lem}

\begin{lem} \label{der_R3_est}
For arbitrary \ $k \geq 2$, \ $\alpha \in (1,2)$  \ and \ $p \in (0,1)$,
 \ there exists \ $C_{6,k}^{\alpha,p} > 0$ \ such that for all \ $n \in \NN$ \ and
 \ $|t| \leq C_2^{\alpha,p} \, n^{1/\alpha}$,
 \[
   \left| \frac{\dd R_{n,3,k}^{\alpha,p}(t)}{\dd t} \right|
   \leq C_{6,k}^{\alpha,p} \frac{|t|^{k-1}}{n^{(k - \alpha)/\alpha}} .
 \]
\end{lem}

\begin{lem} \label{der_R5_7_est}
For arbitrary \ $\ell \in \NN$, \ $\alpha \in (1,2)$ \ and \ $p \in (0,1)$,
 \ there exists \ $C_{7,\ell}^{\alpha,p} > 0$ \ such that for all \ $n \in \NN$
 \ and \ $|t| \leq C_2^{\alpha,p} \, n^{1/\alpha}$,
 \[
   \left| \frac{\dd R_{n,5,\ell}^{\alpha,p}(t)}{\dd t} \right|
   + \left| \frac{\dd R_{n,7,\ell}^{\alpha,p}(t)}{\dd t} \right|
   \leq C_{7,\ell}^{\alpha,p}
        \frac{|t|^{1 + (\ell-1)(2-\alpha)} + |t|^{2\ell - 1}}{n^{\ell(2-\alpha)/\alpha}} .
 \]
\end{lem}

\begin{lem} \label{der_R2_est}
For arbitrary \ $\ell \in \NN$, \ $\alpha \in (1,2)$ \ and \ $p \in (0,1)$,
 \ there exists \ $C_{8,\ell}^{\alpha,p} > 0$ \ such that for all \ $n \in \NN$
 \ and \ $|t| \leq \vare_{\alpha,p} \, n^{1/\alpha}$,
 \[
   \left| \frac{\dd R_{n,2,\ell}^{\alpha,p}(t)}{\dd t} \right|
   \leq C_{8,\ell}^{\alpha,p}
        \frac{|t|^{2\ell-1} \, \ee^{ C_1^{\alpha,p} \, |t|^\alpha / 2}}
             {n^{\ell(2-\alpha)/\alpha}} .
 \]
\end{lem}

\begin{lem} \label{der_R_est}
For arbitrary \ $\ell \in \NN$, \ $\alpha \in (0,1) \cup (1,2)$ \ and
 \ $p \in (0,1)$, \ there exist \ $C_{9,\ell}^{\alpha,p} > 0$ \ such that for all
 \ $n \in \NN$ \ and \ $|t| \leq \vare_{\alpha,p} \, n^{1/\alpha}$,
 \[
   \left| \frac{\dd(\bff_n^{\alpha, p}-\bg_{n,\ell - 1}^{\alpha, p})(t)}{\dd t} \right|
   \leq C_{9,\ell}^{\alpha,p}
        \frac{|t|^{1 + (\ell-1)(2-\alpha)} + |t|^{2\ell - 1 + \alpha}}{n^{\ell(2-\alpha)/\alpha}}
        \ee^{ - C_1^{\alpha,p} \, |t|^\alpha / 2} .
 \]
\end{lem}

By Lemma \ref{der_R_est} and the inequality \eqref{int1}, we obtain
 \ $\tDelta_{n,\ell,3}^{\alpha,p} = O \left( \frac{1}{n^{\ell(2-\alpha)/\alpha}} \right)$.
\ Consequently, we conclude
 \ $\tDelta_{n,\ell}^{\alpha,p}
    = O \left( \frac{1}{n^{\ell(2-\alpha)/\alpha}} + \frac{1}{n^{1 / \alpha}} \right)$,
 \ which implies the first statement.

The reduction of the order \ $O\left(\frac{1}{n^{1/\alpha}}\right)$ \ to 
 \ $o\left(\frac{1}{n^{1/\alpha}}\right)$ \ when \ $\ell \in \{2,3,\dots\}$,
 \ $\alpha \in \left(1, 2 - \frac{1}{\ell} \right) $ \ and
 \ $r^{1/\alpha} \notin \NN$ \ is based again on Lemma \ref{Esseen2}.
By this lemma, there exists a sequence \ $\lambda_n^{\alpha,p}\to\infty$ \ such
 that
 \[
   \int_{\vare_{\alpha,p}}^{\lambda_n^{\alpha,p}}
    \frac{|\bff_{\alpha,p}(s)|^n}{s} \, \dd s
   = o \left( \ee^{-\sqrt{n}/2} \right) \qquad
     \text{as \ $n\to\infty$.}
 \]
Then 
 \begin{align*}
  \tI_{n,1}^{\alpha,p}
  & := \int_{\vare_{\alpha,p} \, n^{1 / \alpha}}^{\lambda_{n-1}^{\alpha,p} \, n^{1/\alpha}}
        \frac{|\bff_n^{\alpha,p}(t)|}{t} \, \dd t
    = \int_{\vare_{\alpha,p}}^{\lambda_{n-1}^{\alpha,p}}
       \frac{|\bff_{\alpha,p}(s)|^n}{s} \, \dd s \\[2mm]
  & \leq \int_{\vare_{\alpha,p}}^{\lambda_{n-1}^{\alpha,p}}
          \frac{|\bff_{\alpha,p}(s)|^{n-1}}{s} \, \dd s
    = o \left( \ee^{-\sqrt{n-1}/2} \right)
    = o \left( \frac{1}{n^{1/\alpha}} \right) \qquad
    \text{as \ $n\to\infty$,}
 \end{align*}
 and similarly,
 \[
   \tI_{n,2}^{\alpha,p}
   := \int_{\vare_{\alpha,p} \, n^{1 / \alpha}}^{\lambda_{n-1}^{\alpha,p} \, n^{1/\alpha}}
       \frac{|\bff_n^{\alpha,p}(t)|}{t^2} \, \dd t
    = o \left( n^{-1/\alpha} \, \ee^{-\sqrt{n-1}/2} \right)
   = o \left( \frac{1}{n^{1/\alpha}} \right) ,
 \]
 as \ $n\to\infty$.
\ By \eqref{f_n},  
 \begin{equation} \label{der_f_n}
   \frac{\dd \bff_n^{\alpha, p}(t)}{\dd t} 
     = n^{(\alpha-1)/\alpha}
       \big(\bff_{\alpha,p}(t/n^{1/\alpha})\big)^{n-1} \,
       \bff_{\alpha,p}'(t/n^{1/\alpha}) \,
       \ee^{- \ii t \mu_1^{\alpha,p} n^{(\alpha-1)/\alpha}} 
       - \ii \, \mu_1^{\alpha,p} \, n^{(\alpha-1)/\alpha} \, \bff_n^{\alpha, p}(t) . 
 \end{equation}
It is easy to check that \ $\bff_{\alpha,p}'(t) = \ii \EE(X \ee^{\ii t X})$ \ for
 all \ $t \in \RR$, \ thus \ $|\bff_{\alpha,p}'(t)| \leq \mu_1^{\alpha,p}$ \ for all
 \ $t \in \RR$, \ hence
 \[
   \tI_{n,3}^{\alpha,p}
   := \int_{\vare_{\alpha,p} n^{1 / \alpha}}^{\lambda_{n-1} n^{1/\alpha}}
       \frac{1}{t} 
       \left| \frac{\dd \bff_n^{\alpha, p}(t)}{\dd t} \right| \, \dd t
    = o \left( n^{(\alpha-1)/\alpha} \, \ee^{-\sqrt{n-1}/2} \right)
   = o \left( \frac{1}{n^{1/\alpha}} \right) ,
 \]
 as \ $n\to\infty$.
\ By Lemma \ref{Osipov} now with \ $T = \lambda_{n-1}^{\alpha,p} \, n^{1/\alpha}$,
 \[
   \tDelta_{n,\ell}^{\alpha,p}
   \leq 2 c_1 \sum_{j=1}^3
               \left( \tDelta_{n,\ell,j}^{\alpha,p} + \tI_{n,j}^{\alpha,p} \right)
        + 2 c_1 \ttI_{n,\ell}^{\alpha,p}
        + c_1 \frac{\tM_{n,\ell}^{\alpha,p}}{\lambda_{n-1}^{\alpha,p} \, n^{1/\alpha}}
   \qquad \text{as \ $n\to\infty$,}
 \]
 where, by the first inequality of Lemma \ref{y}, formula \eqref{der_y}, and
 the statement
 \eqref{int2},
 \begin{align*}
  \ttI_{n,\ell}^{\alpha,p}
  & := \int_{\vare_{\alpha,p} \, n^{1 / \alpha}}^\infty
        \left( \frac{|\bg_{n,\ell-1}^{\alpha,p}(t)|}{t} 
               + \frac{|\bg_{n,\ell-1}^{\alpha,p}(t)|}{t^2}
               + \frac{1}{t} 
                 \left| \frac{\dd \bg_{n,\ell-1}^{\alpha, p}(t)}{\dd t} \right|
        \right) \dd t \\[2mm]
  & = O \left( n^{\ell - 2 - \frac{1}{\alpha}} \,
               \ee^{ - C_1^{\alpha,p} (\vare_{\alpha,p})^\alpha n} \right)
    = o \left( \frac{1}{n^{1/\alpha}} \right) , \qquad
    \text{as \ $n\to\infty$.}
 \end{align*}
Thus \ $\tDelta_{n,\ell}^{\alpha,p} = o \left( \frac{1}{n^{1/\alpha}} \right)$ \ as
 \ $n\to\infty$ \ for
 \ $\alpha \in \left(1, 2 - \frac{1}{\ell} \right)$.

\section{Uniform and nonuniform bounds in asymptotic expansions in local
               merging theorems in the lattice case}

\begin{thm} \label{local_Main_uniform}
For \ $\ell \in \NN$ \ and \ $ r^{ 1 / \alpha} \in \NN $,
 \[
   \sup_{ s \in r^{ 1 / \alpha} \NN }
    \Bigg| \frac{n^{1/\alpha}}{r^{ 1 / \alpha}} \PP \{ S_n = s \}
           - \big( G_{n,\ell - 1}^{\alpha, p} \big)'
             \bigg( \frac{s - c_n^{\alpha,p}}{n^{1/\alpha}} \bigg) \Bigg| 
   = \begin{cases}
      O \left( \frac{1}{n^\ell} \right) ,
       & \text{if \ $ 0 < \alpha < 1 $,} \\[2mm]
      O \left( \frac{[\log_r n]^{2 \ell}}{n^\ell} \right) ,
       & \text{if \ $\alpha = 1$,} \\[2mm]
      O \left( \frac{1}{ n^{ \ell ( 2 - \alpha ) / \alpha } } \right) ,
       & \text{if \ $ 1 < \alpha < 2 $.}
     \end{cases}
 \]
\end{thm}

In case \ $\ell = 1$ \ Theorem \ref{local_Main_uniform} implies for all
 \ $0 < \alpha < 2$ \ the local merging theorem
 \[
   \lim_{n \to \infty}
    \left[ \frac{n^{1/\alpha}}{r^{ 1 / \alpha}} \PP \{ S_n = s \}
           - G_{\alpha, p, \gamma_n}'
             \bigg( \frac{s - c_n^{\alpha,p}}{n^{1/\alpha}} \bigg) \right]
   = 0
 \] 
 for all \ $s \in r^{ 1 / \alpha} \NN$.

\begin{thm} \label{local_MAIN_nonuniform} 
For \ $\ell \in \NN$, \ $\alpha \in (1,2)$ \ and \ $ r^{ 1 / \alpha} \in \NN $,
 \[
   \sup_{ s \in r^{ 1 / \alpha} \NN }
    \bigg( 1 + \frac{|s - \mu_1^{\alpha,p} n|}{n^{1/\alpha}} \bigg)
    \Bigg| \frac{n^{1/\alpha}}{r^{ 1 / \alpha}} \PP \{ S_n = s \}
           - \big( G_{n,\ell - 1}^{\alpha, p} \big)'
             \bigg( \frac{s - \mu_1^{\alpha,p} n}{n^{1/\alpha}} \bigg)
    \Bigg| 
   = O \left( \frac{1}{ n^{ \ell ( 2 - \alpha ) / \alpha } } \right) .
 \]
\end{thm}

\section{Proof Theorems \ref{local_Main_uniform} and
               \ref{local_MAIN_nonuniform}}

First we prove Theorem \ref{local_Main_uniform}.
Fix \ $\ell \in \NN$ \ and \ $\alpha \in (0,2)$ \ such that
 \ $ r^{ 1 / \alpha} \in \NN $.
\ We have
 \[
   \bff_n^{\alpha,p}(t)
   = \EE\big(\ee^{\ii t (S_n - c_n^{\alpha, p}) / n^{1 / \alpha}}\big)
   = \sum_{u \in r^{1 / \alpha} \NN}
      \ee^{\ii t (u - c_n^{\alpha, p}) / n^{1 / \alpha}}
      \PP\{ S_n = u \} ,
 \]
 for all \ $t \in \RR$, \ hence for every \ $s \in r^{1 / \alpha} \NN$,
 \begin{multline*}
  \int_{- \pi n^{1 / \alpha} / r^{1 / \alpha}}^{\pi n^{1 / \alpha} / r^{1 / \alpha}}
   \ee^{- \ii t (s - c_n^{\alpha,p}) / n^{1/\alpha}}
   \bff_n^{\alpha,p}(t) \, \dd t \\[2mm] 
  = \sum_{u \in r^{1 / \alpha} \NN}
     \PP\{ S_n = u \}
     \int_{- \pi n^{1 / \alpha} / r^{1 / \alpha}}^{\pi n^{1 / \alpha} / r^{1 / \alpha}}
      \ee^{\ii t (u - s) / n^{1/\alpha}} \, \dd t
  = \frac{2 \pi n^{1 / \alpha}}{r^{1 / \alpha}} \PP \{ S_n = s \} .
 \end{multline*} 
(In fact, this is the inversion formula for probabilities.)
By Lemma 4 in Cs\"org\H{o} \cite{Cso_07b},
 \[
   \big( G_{n,\ell - 1}^{\alpha, p} \big)'(x)
   = \frac{1}{2 \pi}
     \int_{- \infty}^\infty
      \ee^{- \ii t x} \bg_{n,\ell-1}^{\alpha,p}(t) \, \dd t , \qquad x \in \RR .
 \]
Thus we have
 \begin{align*}
   D_{n,\ell}^{\alpha,p} 
    := \sup_{ s \in r^{ 1 / \alpha} \NN }
         \Bigg| \frac{n^{1/\alpha}}{r^{ 1 / \alpha}} \PP \{ S_n = s \}
                - \big( G_{n,\ell - 1}^{\alpha, p} \big)'
                  \bigg( \frac{s - c_n^{\alpha,p}}{n^{1/\alpha}} \bigg)
         \Bigg| 
    \leq \frac{1}{\pi}
          \left( D_{n,\ell,1}^{\alpha,p}
                 + J_n^{\alpha,p}
                 + J_{n,\ell}^{\alpha,p} \right) ,
 \end{align*}
 where
 \begin{align*}
  D_{n,\ell,1}^{\alpha,p}
  & := \int_0^{\vare_{\alpha,p} \, n^{1 / \alpha}}
        | \bff_n^{\alpha,p}(t) - \bg_{n, \ell - 1}^{\alpha,p}(t) | \, \dd t , \\[2mm] 
  J_n^{\alpha,p}
  & := \int_{\vare_{\alpha,p} \, n^{1 / \alpha}}^{\pi n^{1 / \alpha} / r^{1 / \alpha}}
        | \bff_n^{\alpha,p}(t) | \, \dd t , \\[2mm]
  J_{n,\ell}^{\alpha,p}
  & := \int_{\vare_{\alpha,p} \, n^{1 / \alpha}}^\infty
        | \bg_{n, \ell - 1}^{\alpha,p}(t) | \, \dd t ,
 \end{align*}
 with \ $\vare_{\alpha,p} > 0$ \ from Lemma \ref{R2_est}.

First consider the case \ $\alpha \in (0,1) \cup (1,2)$. 
\ Then by Lemma \ref{R},
 \[
   D_{n,\ell,1}^{\alpha,p}
   = O \left( \frac{1}{n^\ell} + \frac{1}{n^{\ell(2-\alpha)/\alpha}} \right) .
 \]
Moreover,
 \begin{align*}
  J_n^{\alpha,p}
  & = \int_{\vare_{\alpha,p} \, n^{1 / \alpha}}^{\pi n^{1 / \alpha} / r^{1 / \alpha}}
       \left| \bff_{\alpha,p}\left(\frac{t}{n^{1 / \alpha}}\right) \right|^n \dd t
    = n^{1 / \alpha}
      \int_{\vare_{\alpha,p}}^{\pi / r^{1 / \alpha}}
       | \bff_{\alpha,p}(t) |^n \, \dd t \\[2mm]
  & \leq \frac{\pi n^{1 / \alpha}}{r^{1 / \alpha}}
         \left( \sup_{t \in [\vare_{\alpha,p}, \, \pi/r^{1/\alpha}]}
                 |\bff_{\alpha,p}(t)| \right)^n
    = O \left( \frac{1}{n^\ell} + \frac{1}{n^{\ell (2 - \alpha) / \alpha}} \right),
 \end{align*}
 since
 \begin{equation} \label{est_f}
   \sup_{t \in [\vare_{\alpha,p}, \, \pi/r^{1/\alpha}]} |\bff_{\alpha,p}(t)| < 1 .
 \end{equation}
(This follows from the fact that the gain \ $X$ \ in one game takes values also
 in the lattice \ $ r^{ 1 / \alpha} \ZZ$, \ the maximal span of a lattice \ $L$
 \ with \ $\PP \{ X \in L \} = 1$ \ is \ $r^{ 1 / \alpha}$, \ hence
 \ $\sup_{t \in H} |\EE(\ee^{\ii t X})| < 1$ \ for all compact set \ $H \subset \RR$
 \ with $H \cap \frac{2 \pi}{ r^{1/\alpha}} \ZZ = \emptyset$, \ see, e.g., 
 Bhattacharya and Ranga Rao \cite[\S 21]{Bha_Ran_76}.)

By Lemma \ref{y},
 \[
   J_{n,\ell}^{\alpha,p}
   = O \left( n^{\ell - 2} \,
              \ee^{-C_1^{\alpha,p} (\vare_{\alpha,p})^\alpha \, n} \right)
   = O \left( \frac{1}{n^\ell} + \frac{1}{n^{\ell (2 - \alpha) / \alpha}} \right) ,
 \]
 and we obtain the statement for \ $\alpha \in (0,1) \cup (1,2)$.

In case \ $\alpha = 1$, \ in order to estimate \ $D_{n,\ell,1}^{1,p}$, \ we need an
 analogue of Lemma \ref{R}.

\begin{lem} \label{R_1}
For arbitrary \ $\ell \in \NN$ \ and \ $p \in (0,1)$, \ there exist
 \ $C_{5,\ell}^{1,p} > 0$ \ and \ $\vare_{1,p} > 0$, \ such that for all
 \ $n \in \NN$ \ and \ $|t| \leq \vare_{1,p} \, n$,
 \[
   | \bff_n^{1,p}(t) - \bg_{n,\ell - 1}^{1,p}(t) |
   \leq C_{5,\ell}^{1,p}
        \frac{|t|^{\ell+1} + |t|^{2\ell}}{n^\ell}
        \left( 1 + \log_r \frac{2n}{|t|} \right)^{2\ell}
        \ee^{ - C_1^{1,p} |t| / 2} .
 \]
\end{lem}

\begin{proof}
It can be derived in a similar way as Lemma \ref{R}.
First recall that there exists \ $C_1^{1,p} > 0$ \ such that, uniformly in
 \ $\gamma \in (q,1]$, \ we have
 \ $|y_\gamma^{1,p}(t)| \leq C_1^{1,p} (1 + \log_r |t|) \, |t|$ \ for all
 \ $t \in \RR$, \ see Cs\"org\H{o} \cite[Lemma 3]{Cso_07b}.
There exists \ $C_1^{1,p} > 0$ \ such that, uniformly in \ $\gamma \in (q,1]$,
 \ we have \ $\Re ( y_\gamma^{1,p}(t) ) \leq - C_1^{1,p} |t|$ \ for all
 \ $t \in \RR$, \ see inequality (4) in the proof of Lemma 3 of Cs\"org\H{o}
 \cite{Cso_02}.
Next one can show, as in Lemma \ref{R1_x}, that there exist \ $C_2^{1,p} > 0$
 \ and \ $C_3^{1,p} > 0$ \ such that for all \ $n \in \NN$ \ and
 \ $|t| \leq C_2^{1,p} n$,
 \[
   \frac{|x_n^{1,p}(t)|}{n} \leq \frac{1}{2} , \qquad
   |x_n^{1,p}(t)| \leq C_3^{1,p} |t| \left( 1 + \log_r \frac{2n}{|t|} \right) .
 \]
The expression for \ $\bff_n^{1,p}(t)$ \ if \ $|t| \leq C_2^{1,p} n$ \ has the same
 form as in case \ $\alpha \in (0,1) \cup (1,2)$ \ by replacing \ $\alpha$ \ by
 \ 1 \ and \ $\mu_1^{1,p}$ \ by \ $p r \, \log_r n$, \ hence
 \[
   \tx_{n,m}^{1,p}(t)
   := y_{\gamma_n}^{1,p}(t)
      + \ii t p r \, \log_r n
      + \sum_{j=2}^m
         \frac{\mu_j^{1,p}\,(\ii t)^j}
              {j!\,n^{ j - 1 }} .
 \]
Lemma \ref{R3_est} is valid for \ $\alpha = 1$ \ as well.
As in Lemma \ref{R5_7}, there exists \ $C_{2,\ell}^{1,p} > 0$ \ such that for all
 \ $n \in \NN$ \ and \ $|t| \leq C_2^{1,p} n$,
 \[
   |R_{n,5,\ell}^{1,p}(t)| + |R_{n,7,\ell}^{1,p}(t)|
   \leq C_{2,\ell}^{1,p}
        \frac{|t|^{\ell + 1} + |t|^{2 \ell}}{n^\ell}
        \left( 1 + \log_r \frac{2n}{|t|} \right)^{2 \ell} .
 \]
As in Lemma \ref{R2_est}, there exist \ $C_{4,\ell}^{1,p} > 0$ \ and
 \ $\vare_{1,p} \in (0,C_2^{1,p}]$ \ such that for all \ $n \in \NN$ \ and
 \ $|t| \leq \vare_{1,p} \, n$,
 \[
   |R_{n,2,\ell}^{1,p}(t)|
   \leq C_{4,\ell}^{1,p}
        \frac{|t|^{2 \ell}}{n^\ell}
        \left( 1 + \log_r \frac{2n}{|t|} \right)^{2 \ell}
        \ee^{C_1^{1,p} |t| / 2} ,
 \]
 and we obtain the statement of the lemma.
\end{proof}

Now by Lemma \ref{R_1}, we obtain \
 $D_{n,\ell}^{1,p} = O \left( \frac{[\log_r n]^{2 \ell}}{n^\ell} \right)$. 
\ Clearly, \ $J_n^{1,p} = O \left( \frac{[\log_r n]^{2 \ell}}{n^\ell} \right)$ \ can
 be proved as in case \ $\alpha \in (0,1) \cup (1,2)$.
\ Consequently,
 \ $J_{n,\ell}^{1,p} = O \left( \frac{[\log_r n]^{2 \ell}}{n^\ell} \right)$, \ and
 we conclude the statement of Theorem \ref{local_Main_uniform} for
 \ $\alpha = 1$.

In order to prove Theorem \ref{local_MAIN_nonuniform}, first we recall that
 \ $\alpha \in (1,2)$ \ implies \ $\EE(X) < \infty$, \ and hence
 \ $\EE\left| \frac{S_n - \mu_1^{\alpha, p} n}{n^{1 / \alpha}} \right| < \infty$,
 \ and the characteristic function \ $\bff_n^{\alpha,p}$ \ of
 \ $\frac{S_n - \mu_1^{\alpha, p} n}{n^{1 / \alpha}}$ \ is differentiable and
 \[
   \frac{\dd \bff_n^{\alpha,p}(t)}{\dd t}
   = \sum_{u \in r^{1 / \alpha} \NN}
      \frac{\ii (u - \mu_1^{\alpha, p} n)}{n^{1 / \alpha}}
      \ee^{\ii t (u - \mu_1^{\alpha, p} n) / n^{1 / \alpha}}
      \PP\{ S_n = u \} , \qquad t \in \RR .
 \]
Consequently, for every \ $s \in r^{1 / \alpha} \NN$,
 \begin{align*}
  &\int_{- \pi n^{1 / \alpha} / r^{1 / \alpha}}^{\pi n^{1 / \alpha} / r^{1 / \alpha}}
    \ee^{- \ii t (s - \mu_1^{\alpha,p} n) / n^{1/\alpha}}
    \left( \frac{\dd \bff_n^{\alpha,p}(t)}{\dd t} \right) \dd t \\[2mm] 
  &\hspace*{20mm}
   = \sum_{u \in r^{1 / \alpha} \NN}
      \frac{\ii (u - \mu_1^{\alpha, p} n)}{n^{1 / \alpha}}
      \PP\{ S_n = u \}
      \int_{- \pi n^{1 / \alpha} / r^{1 / \alpha}}^{\pi n^{1 / \alpha} / r^{1 / \alpha}}
       \ee^{\ii t (u - s) / n^{1/\alpha}} \, \dd t \\[2mm] 
  &\hspace*{20mm}
   = \frac{\ii (s - \mu_1^{\alpha, p} n)}{n^{1 / \alpha}}
     \frac{2 \pi n^{1 / \alpha}}{r^{1 / \alpha}} \PP \{ S_n = s \} .
 \end{align*} 
By Lemmas 4 and 6 in Cs\"org\H{o} \cite{Cso_07b}, \ $\alpha \in (1,2)$
 \ implies
 \[
   \ii x \big( G_{n,\ell - 1}^{\alpha, p} \big)'(x)
   = \frac{1}{2 \pi}
     \int_{- \infty}^\infty
      \ee^{- \ii t x}
      \left( \frac{\dd \bg_{n,\ell - 1}^{\alpha,p}(t)}{\dd t} \right) \dd t ,
   \qquad x \in \RR .
 \]
Thus we have
 \begin{align*}
   \tD_{n,\ell}^{\alpha,p}
   & := \sup_{ s \in r^{ 1 / \alpha} \NN }
         \frac{|s - \mu_1^{\alpha,p} n|}{n^{1/\alpha}}
         \Bigg| \frac{n^{1/\alpha}}{r^{ 1 / \alpha}} \PP \{ S_n = s \}
                - \big( G_{n,\ell - 1}^{\alpha, p} \big)'
                  \bigg( \frac{s - \mu_1^{\alpha,p} n}{n^{1/\alpha}} \bigg)
         \Bigg| \\[2mm]
   & \leq \frac{1}{\pi}
          \left( \tD_{n,\ell,1}^{\alpha,p}
                 + \tJ_n^{\alpha,p}
                 + \tJ_{n,\ell}^{\alpha,p} \right) ,
 \end{align*}
 where
 \begin{align*}
  \tD_{n,\ell,1}^{\alpha,p}
  & := \int_0^{\vare_{\alpha,p} \, n^{1 / \alpha}}
        \left| \frac{\dd \bff_n^{\alpha,p}(t)}{\dd t}
               - \frac{\dd \bg_{n, \ell - 1}^{\alpha,p}(t)}{\dd t} \right| \dd t ,
       \\[2mm] 
  \tJ_n^{\alpha,p}
  & := \int_{\vare_{\alpha,p} \, n^{1 / \alpha}}^{\pi n^{1 / \alpha} / r^{1 / \alpha}}
        \left| \frac{\dd \bff_n^{\alpha,p}(t)}{\dd t} \right| \dd t , \\[2mm]
  \tJ_{n,\ell}^{\alpha,p}
  & := \int_{\vare_{\alpha,p} \, n^{1 / \alpha}}^\infty
        \left| \frac{\dd \bg_{n, \ell - 1}^{\alpha,p}(t)}{\dd t} \right| \dd t .
 \end{align*}
By Lemma \ref{der_R_est} and inequality \eqref{int1}, 
 \ $\tD_{n,\ell,1}^{\alpha,p} = O \left( \frac{1}{n^{\ell(2-\alpha)/\alpha}} \right)$.
\ Using \eqref{der_f_n} and the inequalities \eqref{est_f} and
 \ $|\bff_{\alpha,p}'(t)| \leq \mu_1^{\alpha,p}$, \ $t \in \RR$, \ we obtain
 \ $\tJ_n^{\alpha,p} = O \left( \frac{1}{n^{\ell (2 - \alpha) / \alpha}} \right)$.
\ By Lemma \ref{der_y_est},
 \ $\tJ_{n,\ell}^{\alpha,p} = O \left( \frac{1}{n^{\ell (2 - \alpha) / \alpha}} \right)$,
\ and by Theorem \ref{local_Main_uniform}, we conclude the statement of Theorem
 \ref{local_MAIN_nonuniform}.

\section*{Acknowledgements}
I am thanksful to Professor S\'andor Cs\"org\H{o} who explained me the
 phenomena of merging in generalized St.~Petersburg games.
This research has been supported by the Hungarian Scientific
 Research Funds under Grant No.\ T048544 and T079128.


\begin{thebibliography}{99}

\bibitem{Bha_Ran_76}
\textsc{Bhattacharya, R. N.} and \textsc{Ranga Rao, R.} (1976).
\textit{Normal Approximation and Asymptotic Expansions}.
Wiley Series in Probability and Mathematical Statistics. 
John Wiley \& Sons, New York-London-Sydney.
MR{0436272}

\bibitem{Bik_66}
\textsc{Bikelis, A.} (1966). 
Estimates of the remainder term in the central limit theorem.
\textit{Litovsk.\ Mat.\ Sb.}
\textbf{6(3)} 323--346.
MR{0210173}

\bibitem{Cri_Wolf_92}
\textsc{Christoph, G.} and \textsc{Wolf, W.} (1992).
\textit{Convergence Theorems with a Stable Limit Law.}
Mathematical Research 70, 
Akademie Verlag, Berlin.
MR{1202035}

\bibitem{Cso_02}
\textsc{Cs\"org\H{o}, S.} (2002). 
Rates of merge in generalized St.~Petersburg games.
\textit{Acta Sci.\ Math.\ (Szeged)}
\textbf{68} 815--847.
MR{1954550}

\bibitem{Cso_03}
\textsc{Cs\"org\H{o}, S.} (2003). 
Merge rates for sums of large gains in generalized St.~Petersburg games.
\textit{Acta Sci.\ Math.\ (Szeged)}
\textbf{69} 441--454.
MR{1992319}

\bibitem{Cso_05}
\textsc{Cs\"org\H{o}, S.} (2005). 
A probabilistic proof of Kruglov's theorem on tails of infinitely divisible
 distributions.
\textit{Acta Sci.\ Math.\ (Szeged)}
\textbf{71} 405--415.
MR{2160375}

\bibitem{Cso_07a}
\textsc{Cs\"org\H{o}, S.} (2007). 
Fourier analysis of semistable distributions.
\textit{Acta Appl.\ Math.}
\textbf{96} 159--175.
MR{2327532}

\bibitem{Cso_07b}
\textsc{Cs\"org\H{o}, S.} (2007). 
Merging asymptotic expansions in generalized St.~Petersburg games.
\textit{Acta Sci.\ Math.\ (Szeged)}
\textbf{73} 297--331.

\bibitem{Cso_Dod_91}
\textsc{Cs\"org\H{o}, S.} and \textsc{Dodunekova, R.} (1991). 
Limit theorems for the Petersburg game.
In: \textit{Sums, Trimmed Sums and Extremes},
Progress in Probability 23, pp.~285--315,
Boston, Birkh\"auser.
MR{1117274}

\bibitem{Esseen_45}
\textsc{Esseen, C.-G.} (1945). 
Fourier analysis of distribution functions.
A mathematical study of the Laplace--Gaussian law.
\textit{Acta Math.}
\textbf{77} 1--125.
MR{0014626}

\bibitem{Hall_82}
\textsc{Hall, P.} (1982).
\textit{Rates of Convergence in the Central Limit Theorem}.
Research Notes in Mathematics 62, 
Pitman, Boston. 
MR{0668197} 

\bibitem{Hall_83}
\textsc{Hall, P.} (1983). 
Fast rates of convergence in the central limit theorem.
\textit{Z. Wahrsch.\ Verw.\ Gebiete}
\textbf{62} 491--507.
MR{0690574}

\bibitem{Katz_68}
\textsc{Katznelson, Y.} (1968). 
\textit{An Introduction to Harmonic Analysis}.
John Wiley \& Sons, Inc., New York-London-Sydney.
MR{0248482}

\bibitem{Kev_09}
\textsc{Kevei, P.} (2009). 
Merging asymptotic expansions for semistable random variables, 
{\it Lithuanian Math.\ J.}
\textbf{49(1)} 40--54.

\bibitem{Mar_85}
\textsc{Martin-L\"of, A.} (1985). 
A limit theorem which clarifies the `Petersburg paradox´.
\textit{J. Appl.\ Probab.}
\textbf{22} 634--643.
MR{0799286}

\bibitem{Osi_67}
\textsc{Osipov, L. V.} (1967). 
Asymptotic expansions in the central limit theorem.
\textit{Vestnik Leningrad.\ Univ.}
\textbf{1967(19)} 45--62.
MR{0216552}

\bibitem{Osi_72}
\textsc{Osipov, L. V.} (1972). 
Asymptotic expansions of the distribution function of a sum of random variables
 with non uniform estimates for the remainder term.
\textit{Vestnik Leningrad.\ Univ.}
\textbf{1972(1)} 51--59.
MR{0300324}

\bibitem{Osi_Pet_67}
\textsc{Osipov, L. V.} and \textsc{Petrov, V. V.} (1967). 
On the estimation of the remainder term in the central limit theorem.
\textit{Teor.\ Verojatnost.\ i Primenen.}
\textbf{12} 322--329.
MR{0216552}

\bibitem{Pet_75}
\textsc{Petrov, V. V.} (1975).
\textit{Sums of Independent Random Variables}.
Ergebnisse der Mathematik und ihrer Grenzgebiete, Band 82. 
Springer-Verlag, New York-Heidelberg. 
MR{0388499}

\end{thebibliography}
\end{document}